\documentclass[draft]{amsart}

\usepackage{graphics}
\usepackage[all]{xy}

\usepackage{amssymb,amsmath}
\usepackage{psfrag}

\usepackage{mathrsfs} 
\usepackage{amsfonts}
\usepackage{amscd}

\usepackage{graphicx}
\usepackage{epstopdf}
\usepackage{tikz}
\usetikzlibrary{arrows,automata}

\usepackage{graphicx,epsfig,float}

\usepackage{url}

\makeatletter
\def\url@leostyle{%
  \@ifundefined{selectfont}{\def\UrlFont{\sf}}{\def\UrlFont{\small\ttfamily}}}
\makeatother
\newtheorem{thm}{Theorem}[section]
\newtheorem{lem}[thm]{Lemma}
\newtheorem{prop}[thm]{Proposition}
\newtheorem{cor}[thm]{Corollary}
\newtheorem{fact}[thm]{Fact}

\newtheorem{clm}[thm]{Claim}

\newtheorem{defn}[thm]{Definition}

\newtheorem{nrmk}[thm]{Remark}
\newtheorem{expl}[thm]{Example}

\newcommand{\pf}{{\bf Proof. }}

\newcommand{\into}{\longrightarrow}

\renewcommand{\tilde}{\widetilde}
\renewcommand{\bar}{\overline}

\newcommand{\NN}{\mathbb{N}}
\newcommand{\ZZ}{\mathbb{Z}}
\newcommand{\QQ}{\mathbb{Q}}

\newcommand{\bA}{{\mathbf A}}

\newcommand{\bI}{{\mathbf I}}
\newcommand{\bJ}{{\mathbf J}}

\newcommand{\maI}{\mathcal I}
\newcommand{\maJ}{\mathcal J}
\newcommand{\maA}{\mathcal A}
\newcommand{\maB}{\mathcal B}
\newcommand{\maC}{\mathcal C}
\newcommand{\maD}{\mathcal D}
\newcommand{\maK}{\mathcal K}
\newcommand{\maY}{\mathcal Y}

\newcommand{\od}{^{\mathrm{op}}}

\newcommand{\oO}{\mbox{${\mathcal O}$}}
\newcommand{\Or}{\mbox{${\mathcal Or}$}}

\renewcommand{\mod}{\mathrm{Mod}}

\newcommand{\cov}{\mathrm{Cov}}

\newcommand{\ob}{\mathrm{Ob}}
\newcommand{\op}{\mathrm{Op}}

\newcommand{\df}{\mathrm{def}}

\newcommand{\Ho}{\mathrm{Hom}}

\newcommand{\id}{\mathrm{id}}

\newcommand{\lind}[1]{\underset{#1}{\underrightarrow{\lim}}}

\begin{document}

\title {On  Pillay's conjecture in the general case}

\author {M\'ario J. Edmundo}

\address{ Departamento de Matem\'atica\\
Faculdade de Ci\^encias,  Universidade de Lisboa\\
Campo Grande, Edif\'icio C6\\ 
P-1749-016 Lisboa, Portugal}

\email{mjedmundo@fc.ul.pt}

\author{Marcello Mamino}

\address{Laboratoire d'Informatique de l'\'Ecole Polytechnique (LIX)\\
B\^atiment Turing, bureau 2011\\
1 Rue Honor\'e d'Estienne d'Orves, Campus de l'\'Ecole Polytechnique\\
91120 Palaiseau, France}

\email{mamino@lix.polytechnique.fr}

\author{Luca Prelli}

\address{ Dipartimento di Matematica\\
Universit\`a degli Studi di Padova\\
via Trieste 63, 35121, Padova, Italia}

\email{lprelli@math.unipd.it}

\author {Janak Ramakrishnan}

\address{Google, Inc. \\
76, 9th Avenue\\
New York NY 10011, USA}

\email{janak@janak.org}

\author{Giuseppina Terzo}

\address{Dipartimento di Matematica e Fisica\\ 
Universit\`a degli Studi della Campania ``Luigi Vanvitelli''\\
Viale Lincoln 5, 81100 Caserta, Italy}

\email{giuseppina.terzo@unina2.it}

\date{\today}
\thanks{The first author was supported by Funda\c{c}\~ao para a Ci\^encia e a Tecnologia, Financiamento Base 2008 - ISFL/1/209.  The third author is a member of the Gruppo Nazionale per l'Analisi Matematica, la Probabilit\`{a} e le loro Applicazioni (GNAMPA) of the Istituto Nazionale di Alta Matematica (INdAM). This work is part of the FCT projects PTDC/MAT/101740/2008 and PTDC/MAT/122844/2010.\newline
 }


\begin{abstract}
Let ${\mathbb M}$ be an arbitrary o-minimal structure. Let $G$ be a definably compact, definably connected,  abelian definable group of dimension $n$.  Here we  compute: (i) the new  intrinsic o-minimal  fundamental group of $G;$ (ii) for each $k>0$, the $k$-torsion subgroups of $G;$  (iii) the  o-minimal cohomology algebra over ${\QQ}$ of $G.$  As  a corollary we obtain a new uniform proof of Pillay's conjecture, an o-minimal analogue of Hilbert's fifth problem, relating definably compact groups to compact real Lie groups, extending the proof already known in o-minimal expansions of ordered fields.
\end{abstract}

\maketitle

\begin{section}{Introduction}\label{section intro}
In this paper we work in an  arbitrary o-minimal structure ${\mathbb M}=(M,<, (c)_{c\in {\mathcal C}}, $ $ (f)_{f\in  {\mathcal F}}, (R)_{R\in {\mathcal R}})$ and are interested in the geometry of definable groups in ${\mathbb M}.$ We refer the reader to \cite{vdd} for basic o-minimality.  O-minimality is the analytic part of model theory and deals with theories of ordered, hence topological, structures satisfying certain tameness properties. It generalizes PL-geometry (\cite{vdd}),  semi-algebraic geometry (\cite{BCR}) and globally sub-analytic geometry (\cite{kr}, also called finitely sub-analytic in \cite{vdd86}) and it is claimed to be the  formalization of  Grothendieck's notion of tame topology (topologie mod\'er\'ee). See \cite{vdd} and \cite{dmi}.

A definable group in an o-minimal structure ${\mathbb M}$ is a group whose underlying set is a definable set and the graph of the group operation is a definable set. The notion of definably compact is the analogue of the notion of semi-algebraically complete and was introduced by Peterzil and Steinhorn in \cite{ps} - see \cite{br3} and the Notes and comments in \cite[page 106]{vdd}. The theory of definable groups, which includes real algebraic groups and semi-algebraic groups, began with Pillay's paper \cite{p1} and has since then grown into a well developed branch of mathematics. The literature contains many interesting results about definable groups which have an analogue in the theory of Lie groups - see \cite{p1}, \cite{PePiSt00a}, \cite{PePiSt00b}, \cite{PePiSt02}, \cite{Ed03}, \cite{eo} and \cite{Ed05}. All these fundamental results hinted at a deeper connection between definably compact definable groups and compact real Lie groups, which was finally formulated in the paper \cite{Pi04} by Pillay. Pillay's conjecture is a non-standard analogue of Hilbert's fifth problem for locally compact topological groups. Roughly it says that after taking the quotient by a ``small subgroup'' (a smallest type-definable subgroup of bounded index)  the quotient when equipped with the  so called logic topology is a compact real Lie group of the same dimension. For more on definable groups and on Pillay's conjecture see \cite{ot} and \cite{Pe07}.

Pillay's conjecture was solved in the following cases: (i) o-minimal expansions of fields \cite{HPePi08} using  new model-theoretic tools and the computation of $m$-torsion subgroups of definably compact abelian groups \cite{eo} (based on o-minimal singular (co)homological arguments); (ii) linear o-minimal expansions of ordered groups using direct methods \cite{ElSt07}; (iii) semi-bounded non-linear o-minimal expansions of ordered groups \cite{Pe09} by reduction to the field case using a refinement of the dichotomy bounded/unbounded for semi-bounded sets studied in \cite{e1}, namely the dichotomy short/long. 

Note that by cases (i), (ii) and (iii) above, Pillay's conjecture was therefore known to be true  for definably compact groups in all o-minimal expansions of ordered groups. What remained open was the case of  o-minimal structures which are not expansions of ordered groups. By   \cite[Theorem 3]{epr}, in the general case we can assume that the definably compact  group lives in a  product of  definable group-intervals. Furthermore, even here, what remained open was to determine the dimension of the quotient of the definably compact group by the ``small subgroup'', which itself reduces to computing the $m$-torsion subgroups of  definably compact abelian groups. 

Here we extend the computation of $m$-torsion subgroups of definably compact abelian groups in o-minimal expansions of real closed fields (\cite{eo}) to arbitrary o-minimal structures. We use o-minimal sheaf cohomology instead of o-minimal singular cohomology and after defining a new o-minimal fundamental group in arbitrary o-minimal structures extending the one from o-minimal expansions of ordered groups we obtain:

\begin{thm}[Structure Theorem]\label{thm main thm}
Let $G$ be a definably compact, definably connected,  abelian definable group of dimension $n$.  Then,  
\begin{enumerate}
\item[(a)] 
the intrinsic o-minimal  fundamental group of $G$ is isomorphic to  ${\ZZ}^n$;

\item[(b)]  
for each $k>0$, the $k$-torsion subgroup of $G$ is isomorphic to  $({\ZZ}/k{\ZZ})^n$, and 
 
\item[(c)]   
the  o-minimal cohomology algebra over ${\QQ}$ of $G$ is isomorphic to the exterior algebra over ${\QQ}$ with $n$ generators of degree one.
\end{enumerate}
\end{thm}

As pointed out in \cite{HPePi08} (see Remark 4 and the end of Section 8), the proof of Pillay's conjecture given in that paper requires the presence of an ambient real closed field only in the following  two places: (i) in the  computation of $m$-torsion subgroups of definably compact abelian groups \cite{eo} and, (ii)  in the following fact on the theory of generic definable subsets first proved in o-minimal expansions of real closed fields in \cite[Theorem 2.1]{PePi07}, using heavily  the work of Dolich \cite{dol} on dividing and forking in o-minimal structures . 

\begin{fact}\label{fact gen}
Let $G$ be a definably compact group defined over a small model ${\mathbb M}_0.$  If $X\subseteq G$ is a closed definable subset,  then the set of ${\mathbb M}_0$-conjugates of $X$ is finitely consistent if and only if $X$ has a point in ${\mathbb M}_0$.
\end{fact}

Since Fact \ref{fact gen} was established in \cite[Theorem 3.2]{et} in arbitrary o-minimal structures after it was generalized to o-minimal expansions of ordered groups in \cite{Pe09} (see point 1 at the beginning of Section 8), we also now have Pillay's conjecture proved in arbitrary o-minimal structures:

\begin{thm}[Pillay's conjecture]\label{thm pi conj}
Let $G$ be a definable group in a $\kappa $-saturated o-minimal structure ${\mathbb M}$ ($\kappa $ large). Then:
\begin{enumerate}
\item
$G$ has a smallest type-definable normal subgroup of bounded index $G^{00}$.
\item
$G/G^{00}$, equipped with the logic topology, is isomorphic, as a topological group, to a compact real Lie group.
\item
If $G$ is definably compact, then $\dim _{{\rm Lie}}(G/G^{00})=\dim _{{\mathbb M}}(G).$\\
\end{enumerate}
\end{thm}

We now explain  the details of the proof of our main result, Theorem \ref{thm main thm}, pointing out to the reader the important points and techniques. 

The strategy is the same as that of the proof of its analogue in o-minimal expansions of real closed fields (\cite{eo}), but we have to use o-minimal sheaf cohomology (\cite{ejp1}) instead of the o-minimal singular homology (\cite{Wo}) and cohomology (\cite{ew}) as well as a new o-minimal fundamental group in arbitrary o-minimal structures generalizing  the one  from o-minimal expansions of fields (\cite{bo}) or ordered groups (\cite{ElSt07}, \cite{EdPa}).

Let $G$ be a definably compact, definably connected, abelian definable  group of dimension  $n.$ 

From the o-minimal (co)homology side we need: (i) the K\"unneth formula to show that the cohomology of $G$ with coefficients in ${\mathbb Q}$ is a graded Hopf algebra of finite type; (ii) the theory of o-minimal (${\mathbb Z}$-)orientability to show that $G$ is orientable and so the (co)homology of $G,$ with coefficients in ${\mathbb Z},$ in degree $n$ is ${\mathbb Z};$ (iii) degree theory for continuous definable maps between orientable definably compact manifolds. These three parts in combination with the fact that the definable homomorphism $p_k:G\to G:x\mapsto kx$ is a definable covering map, gives a lower bound on the size $\#G[k]$ of the subgroup of $k$-torsion points of $G$ of the form $k^r\leq \#G[k]$ where $r$ is the number of generators of the Hopf algebra of $G.$ 

From the o-minimal fundamental group side we need: (iv) the new o-minimal fundamental group is well-connected with the theory of definable covering maps, giving us that $G[k]\simeq ({\mathbb Z}/k{\mathbb Z})^s$ where $s$ is the number of generators of the new o-minimal fundamental group of $G;$ (v) the Hurewicz theorem relating the o-minimal fundamental group  with the o-minimal cohomology in degree one. 

The Hurewicz and the universal coefficients theorems (from the cohomology side) show that $s\leq r$ and so, since we have $k^r\leq k^s,$ we obtain $r=s.$ Since also the sum of the degrees of the $r$ generators of the Hopf algebra of $G$ must be $n,$ because the cohomology of $G$ in degree $n$ is ${\mathbb Q},$ we obtain that $r=s=n$ as required.

Given the above strategy let us now point out exactly which difficulties we had to face in order to implement it.

(i) The K\"unneth formula for the o-minimal singular homology is rather easy from the definitions as in the classical topological case (see \cite{eo} for details). The K\"unneth formula for o-minimal sheaf cohomology (even with coefficients in  constant sheaves) turned out to be rather complicated and is obtained only after the formalism of the Grothendieck six operations on o-minimal sheaves is developed. This formalism was developed in the  recent paper \cite{ep3}, but for definable spaces in   full subcategories $\bA$ of the category of definable spaces such that: \label{page conds a}
\begin{itemize}
\item[(A0)]
cartesian products of objects of $\bA$ are objects of $\bA$ and locally closed definable subsets of objects of $\bA$ are objects of $\bA;$
\item[(A1)]
in every object of $\bA$  every open definable subset is a finite union of open and definably normal definable subsets;
\item[(A2)]
every object of $\bA$ has a definably normal definable completion in $\bA.$
\end{itemize}
Moreover, K\"unneth formula holds for objects $X$ of such a subcategory $\bA$ if furthermore:
\begin{itemize}
\item[(A3)]
for every elementary extension  ${\mathbb S}$ of ${\mathbb M}$ and every sheaf $F$ on the o-minimal site on $X$  we have an isomorphism
$$H^*_c(X; F)\simeq H^*_c(X({\mathbb S}); F({\mathbb S}))$$
where $H_c^*$ is the o-minimal cohomology with definably compact supports (\cite[Example 2.10 and Definition 2.12]{ep1}). 
\end{itemize}
Therefore,  in order to use the K\"unneth formula here, we had to show that:
\begin{itemize}
\item[(*)] the full subcategory of locally closed definable subsets of definably compact definable groups satisfies conditions (A0), (A1) and (A2) and definably compact groups satisfy condition (A3). 
\end{itemize}

\medskip
(ii) O-minimal ${\mathbb Z}$-orientability theory is rather technical both with o-minimal singular homology (\cite{bo}, \cite{beo}) and with o-minimal sheaf cohomology. The difficult part being the proof of the existence of relative fundamental classes associated to orientations. Here this is obtained using a consequence of the  o-minimal Alexander duality theorem proved in \cite{ep3}. 
See Fact \ref{fact alexander dual} of Subsection \ref{subsection prelim rmks} and Definition \ref{defn fund class} of Subsection \ref{subsection degrees}.

\medskip
(iii) Having a good orientation theory available, degree theory is rather classical. The novelty here, see Subsection \ref{subsection degrees},  is that we work with the o-minimal Borel-Moore homology. In any case,   since we only need to work in homology groups of top degree, we actually  don't introduce formally the o-minimal Borel-Moore homology and use instead its description 
 given by the  o-minimal Alexander duality theorem, as the ${\mathbb Z}$-dual of the relative o-minimal cohomology group in top degree.

In both cases, in \cite{eo} and here, the existence of relative fundamental classes associated to orientations  
depends crucially on the existence of finite covers by open definable subsets of  definable manifolds for which we can compute, in \cite{eo} some relative o-minimal singular cohomology groups, and here,  the o-minimal  cohomology with definably compact supports. See Subsection \ref{subsection coho jmcells}. 

Therefore, we had to show that:
\begin{itemize}
\item[(**)]  definably compact definable groups have such finite covers by open definable subsets and, have o-minimal orientation sheaves and  are orientable.
 \end{itemize}

\medskip
(iv),  (v) The existence of an o-minimal fundamental group in arbitrary o-minimal structures extending the o-minimal fundamental group from o-minimal expansions of fields (\cite{bo}) or ordered groups (\cite{ElSt07}, \cite{EdPa}) is one of the main novelties of this paper. See Subsections \ref{subsection new fund grp} and \ref{subsection coverings}.

As observed in the concluding remarks of the paper \cite{eep},  this new o-minimal fundamental group, when relativized to a  full subcategory ${\bf P}$ of the category of locally definable spaces, will have all the properties proved in \cite{eep} (including the good connection to definable covering maps and a Hurewicz theorem) if the following hold:
\begin{itemize}
\item[(P1)]
\begin{itemize}
\item[(a)]
every object of ${\bf P}$ which is definably connected is uniformly definably path connected;
\item[(b)]
definable paths and definable homotopies in objects of ${\bf P}$ can be lifted uniquely to locally definable coverings of such objects;
\end{itemize}
\item[(P2)]
Every object of ${\bf P}$ has  admissible covers by definably simply connected, open definable subsets refining any admissible cover by open definable subsets.
\end{itemize} 

Therefore, we also had to show that:
\begin{itemize}
\item[(***)]  definably compact definable groups live in such subcategories ${\bf P}$ on which the relativization of the new o-minimal fundamental group has  properties (P1) and (P2).
 \end{itemize}

\medskip
The main tool we use to obtain (*), (**) and (***), see Section \ref{section  torsions},  is a consequence of the following result (\cite[Theorem 3]{epr}):

\begin{fact}\label{fact grp epr intro}
If $G$ is a definable group, then there is a definable injection $G\to \Pi _{i=1}^mJ_i$, where each $J_i\subseteq M$ is a definable group-interval.
\end{fact}


In order to take advantage of this fact we develop quite extensively o-minimal topology in the  context of cartesian products of definable group-intervals. This is achieved in most cases by  extending  to the context of cartesian products of definable group-intervals  some techniques used by Berarducci and Fornasiero (\cite{bf}) in o-minimal expansions of ordered groups.  See Section \ref{section top prod grp-int} and Subsection \ref{subsection coho jmcells}. These techniques from \cite{bf} were already used  in \cite[Theorem 1.1]{ep2} to prove  (A3) for definably compact definable groups. \\

\medskip
\noindent
{\bf Acknowledgements:} Various ideas on this paper came already from earlier work by the first author with Gareth Jones and Nicholas Peatfield and also with Giuseppina Terzo - but quite a few details were missing (see for example the preprint referred to in the papers \cite{Pe07} and \cite{Pe09}). The essential work was done under the FCT projects PTDC/MAT/101740/2008 and PTDC/MAT/122844/2010 during a period where we had a very lively model theory group based at CMAF, a research center of the University of Lisbon. We wish to thank all the other members of the Lisbon model theory group at CMAF, Alex Usvyatsov, Pantelis Eleftheriou, Tamara Servi and  Ayhan Gunaydin, for holding reading seminars and discussions on the topic of the paper. We    wish to thank also all of the consultants of the above projects, but especially, Kobi Peterzil, Alessandro Berarducci, Sergei Starchenko,  Anand Pillay, Charlie Steinhorn and Jean-Philippe Rolin for visiting us in Lisbon and talk about this work and other things.

\end{section}

\begin{section}{Preliminaries}\label{section prelimp}
In this section we prove a cover by open cells result and define a new o-minimal fundamental group.

\begin{subsection}{A general cover by open cells result}\label{subsection coverings by open cells}
Here we show that in o-minimal structures with definable choice functions every open  definable subset is a finite union of open definable subsets each definably homeomorphic, by reordering  of coordinates, to an open cell.\\

The following is obtained from the definition of cells (\cite[Chapter 3, \S2]{vdd}):\\

\begin{nrmk}\label{nrmk cells perm}
{\em
Let $C\subseteq M^n$ be a $d$-dimensional cell. Then by definition of cells,  $C$ is a $(i_1, \ldots , i_n)$-cell for some unique sequence $(i_1, \ldots , i_n)$ of $0$'s and $1$'s and there are  $\lambda (1)< \dots <\lambda (d)$  indices $\lambda \in \{1, \ldots , n\}$ for which $i_{\lambda }=1.$ Moreover, if  
$$p_{(i_1, \ldots , i_n)}:M^n\to M^d:(x_1, \ldots , x_n)\mapsto (x_{\lambda (1)}, \ldots , x_{\lambda (d)})$$ 
is the projection, then $C':=p_{(i_1, \ldots , i_n)}(C)$ is an open $d$-dimensional cell in $M^d$ and the restriction $p_C:=p_{(i_1, \ldots , i_n)|C}:C\to C'$ is a definable homeomorphism  (\cite[Chapter 3, (2.7)]{vdd}). 

Let $\tau (1)<\dots <\tau (n-d)$ be the indices $\tau \in \{1, \ldots , n\}$ for which $i_{\tau }=0.$ For each such $\tau $, by definition of cells, there is a definable continuous function $h_{\tau }: \pi _{\tau -1}(C)\subseteq M^{\tau -1}\to M$ where, for each $k=1, \ldots , n$,  $\pi _{k}:M^n\to M^{k}$ is the projection onto the first $k$-coordinates. Moreover we have $\pi _{\tau }(C)=\{(x,h_{\tau }(x)): x\in \pi _{\tau -1}(C)\}.$

Let $f=(f_1, \ldots , f_{n-d}):C'\to M^{n-d}$ be the definable continuous map where for each $l=1, \ldots , n-d$ we set $f_l=h_{\tau (l)}\circ \pi _{\tau (l) -1}\circ p_C^{-1}.$ Let  $\sigma :M^n\to M^n:(x_1,\ldots , x_n)\mapsto (x_{\lambda (1)}, \ldots , x_{\lambda (d)}, x_{\tau (1)}, \ldots , x_{\tau (n-d)}).$ Then we clearly have 
$$\sigma (C)=\left\{\left(x,f(x)\right) : x\in C'\right\}.$$
}
\end{nrmk}

\medskip

\begin{thm}\label{thm open cells}
Suppose that ${\mathbb M}$ has definable choice functions. Let $U$ be an open definable subset of $M^n$. 
Then $U$ is a finite union of open definable sets definably homeomorphic, by reordering of coordinates,  to open cells.
\end{thm}

\pf
It suffices to prove the following. If $C\subset U$ is a cell, then there are finitely many open subsets of~$U$ definably homeomorphic, by reordering of coordinates,  to open cells such that $C$ is contained in the union of them. Since a cell of dimension $n$ is an open set, we assume that the dimension of $C$ is smaller than $n.$

We proceed by induction on the dimension of~$C$. The zero-dimensional case is immediate. Let $C$ be $d$-dimensional and assume the statement for cells of lower dimension. 

Modulo reordering of the coordinates (Remark \ref{nrmk cells perm}) we may assume
\[
C=\left\{\left(x,f(x)\right) : x\in C'\right\}
\]
where $C'\subset M^d$ is a $d$-dimensional open cell and $f\colon C'\subseteq M^d \to M^{n-d}$ is a continuous definable function. Since $U$ is open, for every $x\in C'$ there are $u,v\in M^{n-d}$ such
that $u_i<f_i(x)<v_i$ for every $i=1,\dotsc ,n-d$,
and
\[
\{x\}\times[u_1,v_1]\times\dotsc\times[u_{n-d}, v_{n-d}] \subset U
\]
By definable choice there are definable functions $g=(g_1,\ldots , g_{n-d})\colon C'\to M^{n-d}$ and $ h=(h_1,\ldots , h_{n-d})\colon C'\to M^{n-d}$ such that for every $x\in C'$ we have $g_i(x)<f_i(x)<h_i(x)$ for every $i=1\dotsc n-d$,
and
\[
\{x\}\times[g_1(x),h_1(x)]\times\dotsc\times[g_{n-d}(x), h_{n-d}(x)] \subset U.
\]

Let $O\subset C'$ be the definable set of the continuity points of~$g$ and~$h$. Let $C_1,\dotsc , C_m$ be the $d$-dimensional
cells of a cell decomposition of $C'$ compatible with~$O$. Then for each $i=1,\dotsc ,m$ define
\[
V_i = \left\{ (x, z_1, \ldots , z_{n-d}):  x\in C_i \,\,\textrm{and}\,\, g_l(x)<z_l<h_l(x)\,\,\textrm{for}\,\, l=1, \ldots , n-d\right\}
\]
These sets are clearly definably homeomorphic to open cells (in fact they \textit{are} open cells), 
and $\bigcup_i V_i$ covers $\{(x, f(x)): x \in \bigcup_i C_i\subseteq C'\}\subseteq C.$ Hence we conclude by the induction hypothesis observing that $\dim \left(\{(x, f(x)): x \in C'\setminus \bigcup_i C_i\} \right)=
 \dim \left(C'\setminus \bigcup_i C_i \right) < d$.
\qed \\


\end{subsection}

\begin{subsection}{A general o-minimal fundamental group functor}\label{subsection new fund grp}
Here we introduce an o-minimal fundamental group functor in arbitrary o-minimal structures. We also prove some basic properties of this new general o-minimal fundamental group.\\

First we recall the definition of the category of locally definable manifolds with continuous locally definable maps.\\

A {\it locally definable manifold (of dimension $n$)} is a triple $(S,(U_i,\theta _i)_{i\leq \kappa})$ 
where:

 \begin{itemize}
\item[$\bullet  $]
$S=\bigcup _{i\leq \kappa}U_i$;

\item[$\bullet  $]
each $\theta  _i:U_i\rightarrow  R^{n}$ is an injection such that $\theta _i(U_i)$ is an open definable subset of $M^{n}$;

\item[$\bullet  $]
for all $i, j$, $\theta  _i(U_i\cap U_j)$ is an open definable subset of  $\theta _i(U_i)$ and  the transition maps $\theta _{ij}:\theta  _i(U_i\cap U_j)\rightarrow  \theta  _j(U_i\cap U_j):x\mapsto \theta  _j(\theta  _i^{-1}(x))$ are definable homeomorphisms.
\end{itemize}
We call the $(U_i, \theta _i)$'s the {\it definable charts of $S$}. If $\kappa <\aleph _0 $ then $S$ is a {\it definable manifold}. 

A locally definable manifold $S$ is equipped with the topology such that  a subset $U$ of $S$ is  open  if and only if for each $i$, $\theta  _i(U\cap U_i)$ is an open definable subset of $\theta  _i(U_i)$.

We  say that a subset $A$ of $S$ is \emph{definable} if and only if there is a finite $I_0\subseteq \kappa $ such that $A\subseteq \bigcup _{i\in I_0}U_i$ and for each $i\in I_0$,  $\theta _i(A\cap U_i)$ is a definable subset of $\theta  _i(U_i)$. A subset $B$ of $S$ is \emph{locally definable} if and only if for each $i$,  $B\cap U_i$ is a definable subset of $S$.   We say that  a locally definable manifold $S$ is {\it definably connected} if it is not the disjoint union of two open and closed locally definable subsets. 

If ${\mathcal U}=\{U_{\alpha }\}_{\alpha \in I}$ is a cover of $S$ by open locally definable subsets, we say that ${\mathcal U}$  is  {\it admissible} if for each $i\leq \kappa  $, the cover $\{U_{\alpha }\cap U_i\}_{\alpha \in I}$ of $U_i$ admits a finite subcover. If ${\mathcal V}=\{V_{\beta }\}_{\beta \in J}$ is another cover of $S$ by open locally definable subsets, we say that ${\mathcal V}$ {\it refines} ${\mathcal U}$, denoted by ${\mathcal V}\leq {\mathcal U},$ if there is a map $\epsilon :J\rightarrow  I$ such that $V_{\beta }\subseteq U_{\epsilon (\beta )}$ for all $\beta \in J$.

A map $f:X\to  Y$ between locally definable manifolds with definable charts $(U_i,\theta _i)_{i\leq \kappa _X}$ and $(V_j,\delta _j)_{j\leq \kappa _Y}$ respectively is {\it a locally definable map} if for every finite $I\subseteq \kappa _X$ there is a finite $J\subseteq \kappa _Y$ such that:
\begin{itemize}
\item[$\bullet$]  $f(\bigcup _{i\in I}U_i)\subseteq \bigcup _{j\in J} V_j$;

\item[$\bullet$] the  restriction $f_{|}:\bigcup _{i\in I}U_i \to  \bigcup _{j\in J} V_j$ is a definable map between definable manifolds, i.e., for each $i\in I$ and every $j\in J$,
    $\delta _j\circ f\circ \theta _i^{-1}:\theta _i(U_i)\to  \delta _j(V_j)$ is a definable map between definable sets.
\end{itemize}
Thus we have the category of locally definable manifolds with locally definable continuous maps.\\

\begin{defn}
{\em
By a {\it basic $d$-interval}, short for {\it basic directed interval},  we mean a tuple 
$$\maI=\langle [a,b], \langle 0_{\maI}, 1_{\maI}\rangle \rangle$$
where $a,b \in M$ with $a<b$ and $\langle  0_{\maI}, 1_{\maI}\rangle \in \{\langle a, b\rangle , \langle b, a\rangle \}.$ The {\it domain} of $\maI $ is $[a,b]$ and the {\it direction} of $\maI $ is $\langle 0_{\maI}, 1_{\maI}\rangle.$ The {\it opposite} of  $\maI $ is the basic $d$-interval
$$\maI \od =\langle [a,b], \langle 0_{\maI \od}, 1_{\maI \od} \rangle \rangle $$
with the same domain and opposite direction $\langle 0_{\maI \od}, 1_{\maI \od}\rangle =\langle 1_{\maI }, 0_{\maI} \rangle .$

If $\maI _i=\langle [a_i, b_i], \langle 0_{\maI _i}, 1_{\maI _i}\rangle \rangle$ are basic $d$-intervals, for $i=1, \ldots , n,$ we define the {\it $d$-interval}, short for {\it directed interval}, $ \maI _1\wedge \dots \wedge \maI _n,$  whose {\it domain } is the set 
\[
[a_1, b_1]\wedge \dots \wedge [a_n, b_n]:= \raisebox{1ex}{$\mathsurround=0pt\displaystyle \bigsqcup_ i\{c_i\}\times [a_i, b_i]$}
\Big/ \raisebox{-1ex}{$\mathsurround=0pt\displaystyle\sim$}
\]  
where $c_1,\dots , c_n$ are $n$ distinct points of~$M$ and $\sim $ is the equivalence relation defined by  $(c_i,1_{\maI _i})\sim(c_{i+1}, 0_{\maI _{i+1}})$ for each~$i=1,\dots , n-1$ and identity elsewhere.   The  {\it direction} of $\maI _1\wedge \dots \wedge \maI _n $ is $\langle 0_{\maI _1\wedge \dots \wedge \maI _n}, 1_{\maI _1\wedge \dots \wedge \maI _n} \rangle $ where 
 $0_{\maI _1\wedge \dots \wedge \maI _n}=\langle c_1, 0_{\maI _1} \rangle $ and $1_{\maI _1\wedge \dots \wedge \maI _n}=\langle c_n, 1_{\maI _n}\rangle .$

The {\it opposite} of  $\maI _1\wedge \dots \wedge \maI _n $ is the $d$-interval
$$(\maI _1\wedge \dots \wedge \maI _n ) \od =\langle [a_1, b_1]\wedge \dots \wedge [a_n, b_n], \langle 0_{(\maI _1\wedge \dots \wedge \maI _n) \od}, 1_{(\maI _1\wedge \dots \wedge \maI _n) \od} \rangle \rangle $$
with the same domain and opposite direction 
\[
\langle 0_{(\maI _1\wedge \dots \wedge \maI _n) \od}, 1_{(\maI _1\wedge \dots \wedge \maI _n) \od}\rangle =\langle 1_{\maI _1\wedge \dots \wedge \maI _n }, 0_{ \maI _1\wedge \dots \wedge \maI _n} \rangle .
\]
}
\end{defn}

\begin{fact}\label{fact op}
{\em
If $\maI _i=\langle [a_i, b_i], \langle 0_{\maI _i}, 1_{\maI _i}\rangle \rangle$ are basic $d$-intervals, for $i=1, \ldots , n,$ then $(\maI _1\wedge \dots \wedge \maI _n ) \od = \maI _n \od\wedge \dots \wedge \maI _1 \od .$
}
\end{fact}

Below, for the notion of {\it definable space} we refer the reader to \cite[page 156]{vdd}.

\begin{lem}\label{lem d-int space}
Let $\maI =\langle I, \langle 0_{\maI}, 1_{\maI}\rangle \rangle $ be a $d$-interval. Then the domain $I$ of $\maI$ is a Hausdorff, definably compact, definable space  of dimension one which is equipped with a definable total order $<_{\maI}.$
\end{lem}

\pf
Let $\maI _i=\langle [a_i, b_i], \langle 0_{\maI _i}, 1_{\maI _i}\rangle \rangle$ be basic $d$-intervals, for $i=1, \ldots , n,$ and suppose that $\maI =\maI _1\wedge \dots \wedge \maI _n.$ Then $I=[a_1, b_1]\wedge \dots \wedge [a_n, b_n]$ is clearly a Hausdorff, definably compact, definable space of dimension one. 

For each $i$ let $<_{\maI _i}$ be the total order on $[a_i, b_i]$ which is $<$ if $\langle 0_{\maI _i}, 1_{\maI _i}\rangle =\langle a_i, b_i\rangle $ or $>$ if $\langle 0_{\maI _i}, 1_{\maI _i}\rangle =\langle b_i, a_i\rangle .$ Then  total ordering on~$\maI$ is given by
$x<_{\maI}y$ if $x\nsim y$ and either $x,y\in [a_i, b_i]$ for some $i$ and~$x< _{\maI _i}y$, or $x\in I_i$ and~$y\in I_j$ with~$i<j$. 
\qed \\


Due to Lemma \ref{lem d-int space}, below we will identify a $d$-interval $\maI =\langle I, \langle 0_{\maI}, 1_{\maI}\rangle \rangle $ with its domain equipped with the definable total order $<_{\maI}.$ In particular, since the domain $I$ of $\maI \od$ is a definable space of dimension one which is equipped with the definable total order $>_{\maI},$ we  have an
order reversing definable homeomorphism (with respect to the topologies given by the orders)
$$o_{\maI}: \maI \to \maI \od$$
given by the identity on the domain. \\


Given two $d$-intervals~$\maI=\maI _1\wedge\dotsb\wedge \maI _n$ and $\maJ =\maJ _1\wedge\dotsb\wedge \maJ _m$, we define the $d$-interval
\[
\maI\wedge\maJ  = \maI _1\wedge\dotsb\wedge \maI _n\wedge \maJ _1\wedge\dotsb\wedge \maJ _m
\]
and we will regard
$\maI$ and~$\maJ $ as definable subsets of~$\maI \wedge \maJ .$ 

We say that $\maI$ and~$\maJ $ are equal, denoted~$\maI=\maJ $, if $n=m$ and $\maI _i=\maJ _i$ for all~$i=1,\dots , n.$\\

Below, if $X$ be a locally definable manifold and $Y$ is a definable space, we say that $h:Y\to X$ is a {\it definable continuous map} if for some (equivalently, for every) definable subspace $U$ of $X$ with $h(Y)\subseteq U,$ the map $h:Y\to U$ is a definable continuous map between definable spaces.\\

\begin{defn}
{\em
Let $X$ be a locally definable manifold. A {\it  definable path} $\alpha :\maI \to  X$ is a continuous  (with respect to the topology on $\maI$ given by the order) definable map from some $d$-interval $\maI$ to~$X$. We define $\alpha _0:= \alpha (0_{\maI})$ and $\alpha _1:= \alpha (1_{\maI})$ and call the them the end points of the definable path $\alpha .$

A definable path $\alpha :\maI \to  X$ is  {\it constant} if  $\alpha _0=\alpha (t)$ for all $t\in \maI$. Below, given a $d$-interval $\maI$ and a point $x\in X,$ we denote by ${\mathbf c}^x_{\maI}$ the constant definable path in $X$ with endpoints $x.$

A definable path $\alpha :\maI \to  X$ is a {\it definable loop} if  $\alpha _0=\alpha _1.$ The {\it inverse} $\alpha ^{-1}$ of a definable path $\alpha :\maI \to  X$ is the definable path  
$$\alpha ^{-1}:=\alpha \circ o_{\maI}^{-1}:\maI \od \to  X.$$ 
A {\it concatenation of two definable paths} $\gamma :\maI \to  X$  and $\delta :\maJ\to  X$ with $\gamma (1_{\maI})=\delta (0_{\maJ})$ is the  definable path $\gamma \cdot \delta :\maI \wedge \maJ \to   X$ with:
\begin{displaymath}
(\gamma \cdot \delta )(t) =\left \{ \begin{array}{ll}
\gamma (t)
& \textrm{if $t\in \maI$}\\
\,\,\, \\
\delta (t)
& \textrm{if $t\in \maJ$.}
\end{array} \right.
\end{displaymath}

We say that $X$ is {\it definably path connected} if for every $u,v$ in $X$ there is a definable path $\alpha :\maI \to   X$ such that $\alpha _0=u$ and $\alpha _1=v.$ \\
}
\end{defn}

In the special case required for our applications we shall prove later, see Corollary \ref{cor pi cover  def jman} (1), that  being definably connected is equivalent to  being definably path connected.\\

\medskip
Let $X$ be a locally definable manifold and $Y$ a definable space.  Given two definable continuous maps $f,g:Y\to  X$, we say that a definable continuous map $F(t,s):Y\times \maJ\to   X$ is a {\it  definable homotopy between $f$ and $g$} if $f=F_0:=F_{0_{\maJ}}$ and $g=F_1:=F_{1_{\maJ}}$, where $F_s:=F(\cdot ,s)$  for all $s\in \maJ.$ In this situation we say that $f$ and $g$ are {\it definably homotopic}, denoted $f\sim g$.\\

Since definable paths need not have the same domain, the notion of homotopic definable paths is not contained in the notion of homotopic definable maps just defined:

\begin{defn}
{\em
Two definable paths $\gamma :\maI \to   X$, $\delta :\maJ\to   X$, with $\gamma _0=\delta _0 $ and $\gamma _1=\delta _1$, are called {\it definably homotopic}, denoted $\gamma \approx \delta ,$ if there  are $d$-intervals $\maI'$ and~$\maJ '$
such that $\maJ '\wedge\maI=\maJ \wedge\maI'$, and there is a definable homotopy 
$${\mathbf c}_{\maJ '}^{\gamma_0}\cdot \gamma
\sim
\delta\cdot  {\mathbf c}_{\maI'}^{\delta_1}$$
fixing the end points (i.e., they are definably homotopic by a definable homotopy $F:\maK \times \maA \to X,$ where $\maK = \maJ '\wedge\maI=\maJ \wedge\maI', $ such that $F(0_{\maK}, s)=\gamma _0=\delta _0$ and $F(1_{\maK}, s)=\gamma _1=\delta _1$ for all $s\in \maA.$)
}
\end{defn}

\medskip
The goal now is to show that definable homotopy of definable paths $\approx $ is an equivalence relation compatible with concatenation. The next two observations show that definable homotopy $\sim $ is an equivalence relation compatible with concatenation. However we have to do more since the relation $\approx $ does not assume that the domains of the definable paths are the same.\\

\begin{nrmk}\label{th-homotopy-eq}
{\em
Let $X$ be a locally definable manifold and $Y$ a definable space.  Then definable homotopy of  definable continuous maps $Y\to  X$ is an equivalence relation. 

Indeed, $F:Y\times \maJ\to   X: (t,s)\mapsto f(t)$ is a definable homotopy between $f$ and $f$; if  $F:Y\times \maJ\to   X$ is a definable homotopy between $f$ and $g$, then $H:=F\circ ({\rm  id}_Y\times   o^{-1}_{\maJ }):Y\times \maJ \od \to X$ is a definable homotopy between $g$ and $f$; if  $F:Y\times \maJ\to   X$ is a definable homotopy between $f$ and $g$ and if  $G:Y\times \maK \to   X$ is a definable homotopy between $g$ and $h$, then  $H:Y\times (\maJ \wedge \maK )\to   X$ with
\begin{displaymath}
H(t,s) =\left \{ \begin{array}{ll}
F(t,s)
& \textrm{if $s\in \maJ$}\\
\,\,\, \\
G(t,s)
& \textrm{if $s\in \maK$.}
\end{array} \right.
\end{displaymath}
is a definable homotopy between $f$ and $h.$ \\
}
\end{nrmk}

\begin{nrmk}\label{nrmk hom conc same}
{\em
Let $X$ be a locally definable manifold. If $\gamma _i:\maI \to X$ ($i=1,2$) and $\delta  :\maJ \to X$  are  definable paths with $\gamma_1\sim \gamma _2 $  and   $(\gamma _i)_1=\delta _0$ for $i=1,2,$ then
$\gamma_1\cdot \delta  \sim \gamma  _2\cdot \delta .$ 

Let  $F:\maI \times \maA\to   X$ be a definable homotopy between $\gamma _1$ and $\gamma _2.$ Let $i:\maI \to \maI \wedge \maJ$ and $j:\maJ \to \maI \wedge \maJ $ be the obvious definable immersions.
Then  $H:(\maI \wedge \maJ )\times \maA \to   X$ with
\begin{displaymath}
H(t,s) =\left \{ \begin{array}{ll}
F(t',s)
& \textrm{for $t=i(t')$ and $s\in \maA$}\\
\,\,\, \\
\delta (t')
& \textrm{for $t=j(t')$}\\
\end{array} \right.
\end{displaymath}
is a definable homotopy between $\gamma _1\cdot \delta $ and $\gamma _2\cdot \delta .$ 

Similarly, if $ \lambda :\maJ \to X$ is a  definable path with  $\lambda _1=( \gamma _i)_0$  for $i=1,2,$ then
$\lambda \cdot \gamma _1\sim \lambda \cdot \gamma _2.$ 

Therefore, by transitivity of $\sim $ (Remark \ref{th-homotopy-eq}),  if $\delta _i:\maJ \to X$ ($i=1,2$)  are  definable paths with $\delta _1\sim \delta  _2 $  and   $(\gamma _i)_1=(\delta _i)_0$ for $i=1,2,$ then
$\gamma_1\cdot \delta _1 \sim \gamma  _2\cdot \delta _2.$ \\
}
\end{nrmk}

\begin{nrmk}\label{th-constant-shifting}
{\em
Let $X$ be a locally definable manifold. If $\gamma : \maI \to X$ is a definable path and $\maJ$ is any $d$-interval, then
\[
{\mathbf c}_{\maI\wedge\maJ }^{\gamma_0}\cdot \gamma
\sim
\gamma\cdot {\mathbf c}_{\maJ \wedge\maI}^{\gamma_1}.
\]
Indeed let
\begin{align*}
i_1: \maI &\to \maI\wedge(\maJ \wedge\maI) \\
i_2: \maI &\to (\maI\wedge\maJ )\wedge\maI \\
\end{align*}
be the two immersions  of~$\maI$ in~$\maI\wedge\maJ \wedge\maI$.  Then $H:(\maI \wedge \maJ \wedge \maI )\times \maI\to   X$ with
\begin{displaymath}
H(t,s) =\left \{ \begin{array}{ll}
\gamma(t') &\text{for $t=i_1(t')$ with~$t'<s$} \\
\,\,\, \\
\gamma(t') &\text{for $t=i_2(t')$ with~$s<t'$} \\
\,\,\, \\
\gamma(s)  &\text{otherwise}
\end{array} \right.
\end{displaymath}
is the required definable homotopy. \\

}
\end{nrmk}

It follows from Remarks \ref{th-constant-shifting} and \ref{nrmk hom conc same} that:

\begin{nrmk}\label{nrmk sim and approx}
{\em
Let $X$ be a locally definable manifold. If $\delta _i:\maJ \to X$ ($i=1,2$)  are  definable paths with $\delta _1\sim \delta  _2 ,$ then $\delta _1 \approx \delta _2.$ \\
}
\end{nrmk}

\begin{lem}\label{th-weak-homotopy-sym}
Let $X$ be a locally definable manifold. Let   $\gamma :\maI\to X$ and $\delta :\maJ \to X$ be definable paths with $\gamma _0=\delta _0$ and $\gamma _1=\delta _1.$ Then  the following are equivalent:
\begin{enumerate}
\item
$\gamma \approx \delta .$ 
\item
There are four $d$-intervals $\maA$, $\maB$, $\maC $, and~$\maD$, such that $\maA \wedge\maI\wedge\maB  = \maC \wedge\maJ \wedge\maD $ and
$
{\mathbf c}_{\maA }^{\gamma_0}\cdot \gamma \cdot {\mathbf c}_{\maB }^{\gamma_1}
\sim
{\mathbf c}_{\maC }^{\delta_0}\cdot \delta \cdot {\mathbf c}_{\maD }^{\delta_1}.
$
\end{enumerate}
\end{lem}

\pf
Assume (1). Consider $d$-intervals $\maI'$ and~$\maJ '$ such that $\maJ '\wedge\maI=\maJ \wedge\maI'$ and there is a definable homotopy ${\mathbf c}_{\maJ '}^{\gamma_0}\cdot \gamma
\sim \delta\cdot  {\mathbf c}_{\maI'}^{\delta_1}.$ Let $\maA =\maJ ',$ $\maB= \maJ \wedge \maI ',$ $\maC= \maJ '\wedge \maI$ and $\maD =\maI '.$ Then $\maA \wedge\maI\wedge\maB  = \maC \wedge\maJ \wedge\maD $ and 
we have by Remark \ref{nrmk hom conc same} 
\begin{align*}
{\mathbf c}_{\maA }^{\gamma_0}\cdot \gamma \cdot {\mathbf c}_{\maB }^{\gamma_1} &= ({\mathbf c}_{\maJ '}^ {\gamma_0} \cdot  \gamma )\cdot  {\mathbf c}_{\maB }^{\gamma_1} \sim (\delta \cdot {\mathbf c}_{\maI '}^{\delta _1}) \cdot {\mathbf c}_{\maB}^{\delta _1} = (\delta \cdot {\mathbf c}_{\maI '}^{\delta _1}) \cdot {\mathbf c}_{\maJ \wedge \maI'}^{\delta _1} = (\delta \cdot {\mathbf c}_{\maI ' \wedge \maJ}^{\delta _1}) \cdot {\mathbf c}_{\maI'}^{\delta _1}.
\end{align*}
Since
\begin{align*}
{\mathbf c}_{\maC }^{\delta_0}\cdot \delta \cdot {\mathbf c}_{\maD }^{\delta_1} &= ({\mathbf c}_{\maC}^{\delta_0} \cdot  \delta )\cdot   {\mathbf c}_{\maI '}^{\delta_1} = ({\mathbf c}_{\maJ' \wedge \maI}^{\delta_0} \cdot  \delta )\cdot   {\mathbf c}_{\maI '}^{\delta_1} = ({\mathbf c}_{\maJ \wedge \maI '}^{\delta_0} \cdot  \delta )\cdot   {\mathbf c}_{\maI '}^{\delta_1}. 
\end{align*}
We conclude by Remarks~\ref{th-constant-shifting} and again \ref{nrmk hom conc same} and transitivity of $\sim $ (Remark~\ref{th-homotopy-eq}).

Assume (2). Consider four $d$-intervals $\maA$, $\maB$, $\maC $, and~$\maD$, such that $\maA \wedge\maI\wedge\maB  = \maC \wedge\maJ \wedge\maD $ and
$
{\mathbf c}_{\maA }^{\gamma_0}\cdot \gamma \cdot {\mathbf c}_{\maB }^{\gamma_1}
\sim
{\mathbf c}_{\maC }^{\delta_0}\cdot \delta \cdot {\mathbf c}_{\maD }^{\delta_1}.
$
Let 
\begin{align*}
\maJ '&=\maJ \wedge\maA \wedge\maI\wedge\maB  &
\maI'&=\maC \wedge\maJ \wedge\maD \wedge\maI .&
\end{align*}
Then $\maJ '\wedge\maI=\maJ \wedge\maI'$  and by Remark~\ref{th-constant-shifting} we also have
\begin{align*}
{\mathbf c}_{\maJ '}^{\gamma_0}\cdot \gamma & = {\mathbf c}_{\maJ }^ {\gamma_0} \cdot {\mathbf c}_{\maA }^{\gamma_0} \cdot  {\mathbf c}_{\maI}^{\gamma _0} \cdot  {\mathbf c}_{\maB }^{\gamma_0}\cdot \gamma   \sim {\mathbf c}_{\maJ }^ {\gamma_0} \cdot  ({\mathbf c}_{\maA }^{\gamma_0} \cdot  \gamma \cdot  {\mathbf c}_{\maB }^{\gamma_1})
\cdot  {\mathbf c}_{\maI}^{\gamma_1} \\
\delta \cdot  {\mathbf c}_{\maI'}^{\delta_1} & = \delta \cdot    {\mathbf c}_{\maC }^{\delta_1} \cdot  {\mathbf c}_{\maJ }^{\delta_1} \cdot    {\mathbf c}_{\maD }^{\delta_1}\cdot  {\mathbf c}_{\maI }^{\delta_1} \sim {\mathbf c}_{\maJ }^{\delta_0} \cdot 
({\mathbf c}_{\maC }^{\delta_0} \cdot  \delta \cdot  {\mathbf c}_{\maD }^{\delta_1})\cdot  {\mathbf c}_{\maI}^{\delta_1}.
\end{align*}
We conclude by Remark \ref{nrmk hom conc same} and transitivity of $\sim $ (Remark~\ref{th-homotopy-eq}).
\qed \\

\begin{prop}\label{th-weak-homotopy-eq}
Let $X$ be a locally definable manifold and $x_0, x_1\in X$. Let ${\mathbb P}(X, x_0, x_1)$ denote the set of all  definable paths in $X$ that start at $x_0$ and end at $x_1.$ Then  the restriction of  $\approx ,$  the relation of being definably homotopic, to ${\mathbb P}(X, x_0, x_1)\times {\mathbb P}(X, x_0, x_1)$ is an equivalence relation on ${\mathbb P}(X, x_0, x_1)$.
\end{prop}

\pf
For reflexivity, let $\gamma :\maI \to X$ be a definable path in ${\mathbb P}(X, x_0, x_1),$ and take $\maI '=\maI \wedge \maI =\maJ '.$ Then $\maJ '\wedge \maI = \maI \wedge \maI '$ and ${\mathbf c}_{\maJ '}^{\gamma_0}\cdot \gamma
\sim
\gamma \cdot  {\mathbf c}_{\maI'}^{\gamma _1}$ by Remark~\ref{th-constant-shifting}. 

Symmetry follows at once from Lemma~\ref{th-weak-homotopy-sym}. 

For transitivity consider definable paths $\gamma : \maI \to X,$ $\lambda :\maY \to X$ and $\delta :\maJ \to X$ in ${\mathbb P}(X, x_0, x_1)$ and assume that 
$\gamma\approx\lambda$ and $\lambda \approx\delta .$ Then there are $d$-intervals~$\maJ '$, $\maY '$, and~$\maY ''$, $\maI'$, 
such that $\maY '\wedge \maI=\maY \wedge \maI ',$ $\maJ'\wedge \maY =\maJ \wedge \maY ''$ and 
\begin{align*}
{\mathbf c}_{\maY '}^{\gamma_0}\cdot \gamma &\sim \lambda\cdot {\mathbf c}_{\maI '}^{\lambda_1} \\
{\mathbf c}_{\maJ '}^{\lambda_0}\cdot \lambda &\sim \delta \cdot {\mathbf c}_{\maY ''}^{\delta_1}.
\end{align*}
By Remark~\ref{nrmk hom conc same} we have
\begin{align*}
{\mathbf c}_{\maJ '\wedge\maY '}^{\gamma_0}\cdot \gamma &= {\mathbf c}_{\maJ '}^{\gamma _0}\cdot {\mathbf c}_{\maY '}^{\gamma_0}\cdot \gamma  \sim   {\mathbf c}_{\maJ '}^{\gamma _0}\cdot \lambda\cdot {\mathbf c}_{\maI '}^{\lambda_1}\\
 \delta \cdot {\mathbf c}_{\maY ''\wedge\maI '}^{\delta_1} &=  \delta \cdot {\mathbf c}_{\maY ''}^{\delta _1} \cdot {\mathbf c}_{\maI '}^{\delta_1} \sim {\mathbf c}_{\maJ '}^{\lambda_0}\cdot \lambda  \cdot {\mathbf c}_{\maI '}^{\delta_1}.
\end{align*}
Since $\maJ \wedge (\maY ''\wedge \maI')=(\maJ '\wedge \maY ')\wedge \maI$ we conclude that $\gamma \approx \delta.$
\qed \\

\begin{lem}\label{th-weak-homotopy-composition}
Let $X$ be a locally definable manifold. Let $\gamma, \gamma', \delta$ and $\delta'$ be definable paths  in  $X$ such that $\gamma _1=\delta _0$ and $\gamma '_1=\delta '_0.$
If $\gamma\approx\gamma'$ and $\delta\approx\delta'$, then $\gamma \cdot \delta \approx \gamma' \cdot  \delta' .$
\end{lem}

\pf
By transitivity of $\approx $ (Proposition \ref{th-weak-homotopy-eq}) it suffices to prove the case $\delta=\delta'$. Suppose that $\gamma :\maI \to X$, $\gamma' :\maJ \to X$ and~$\delta :\maY \to X.$
By hypothesis there are $d$-intervals $\maJ '$ and~$\maI'$ such that $\maJ '\wedge \maI = \maJ \wedge \maI '$ and ${\mathbf c}_{\maJ '}^{\gamma_0}\cdot \gamma \sim \gamma'\cdot {\mathbf c}_{\maI'}^{\gamma'_1}.$
By Remarks ~\ref{nrmk hom conc same} and \ref{th-constant-shifting} we obtain
\[
{\mathbf c}_{\maJ '}^{\gamma_0}\cdot \gamma \cdot \delta \cdot {\mathbf c}_{\maJ \wedge\maY }^{\delta_1}
\sim
\gamma' \cdot {\mathbf c}_{\maI'}^{\gamma'_1}\cdot \delta \cdot {\mathbf c}_{\maJ \wedge\maY }^{\delta_1}
\sim
{\mathbf c}_{\maJ \wedge\maI'\wedge\maY }^{\gamma'_0} \cdot \gamma' \cdot \delta
\]
and we conclude by Lemma~\ref{th-weak-homotopy-sym} (using also Remark ~\ref{nrmk hom conc same}).
\qed \\

\begin{lem}\label{th-homotopy-inverse}
Let $X$ be a locally definable manifold and let $\gamma : \maI \to X$ be a definable path in $X.$ Then $\gamma \cdot \gamma^{-1}\sim {\mathbf c}_{\maI \wedge \maI \od}^{\gamma _0}$ and so  $\gamma \cdot \gamma^{-1}\approx {\mathbf c}_{\maI \wedge \maI \od}^{\gamma _0}.$ 
\end{lem}

\pf
We have that $H:(\maI \wedge \maI \od)\times \maI \od \to   X$ with  
\begin{displaymath}
H(t,s) =\left \{ \begin{array}{ll}
\gamma(t) & \text{if $t\in\maI$ and $s<o_{\maI}(t)$} \\
\,\,\, \\
\gamma^{-1}(t) & \text{if $t\in \maI \od$ and $s<t$} \\
\,\,\, \\
\gamma^{-1}(s) &\text{otherwise}
\end{array} \right.
\end{displaymath}
is the  definable homotopy $\gamma \cdot \gamma^{-1}\sim {\mathbf c}_{\maI \wedge \maI \od}^{\gamma _0}$ and the rest follows from Remark \ref{nrmk sim and approx}. 
\qed \\

Let $X$ be a locally definable manifold  and $e_X\in X.$ If ${\mathbb L}(X, e_X)$ denotes the set of all definable loops that start and end at a fixed  element $e_X$ of $X$ (i.e. ${\mathbb L}(X,e_X)={\mathbb P}(X,e_X,e_X)$), the restriction of $\approx  $ to ${\mathbb L}(X, e_X)\times {\mathbb L}(X, e_X)$ is an equivalence relation on ${\mathbb L}(X, e_X)$. 

\begin{defn}
{\em 
We define the {\it o-minimal fundamental group} $\pi _1(X, e_X)$ of $X$ by
$$\pi _1(X, e_X):=
\raisebox{1ex}{$\mathsurround=0pt\displaystyle {\mathbb L}(X,e _X)$}
\Big/ \raisebox{-1ex}{$\mathsurround=0pt\displaystyle \approx $}
$$
and we set $[\gamma ]:=$ the class of $\gamma \in {\mathbb L}(X, e_X)$. By Lemmas~\ref{th-weak-homotopy-composition} and ~\ref{th-homotopy-inverse},  $\pi _1(X, e_X)$ is indeed a group with group operation given by $[\gamma ][\delta ]=[\gamma \cdot \delta ]$ and identity the class a of constant loop at $e_X.$ 

If $f:X\to  Y$ is a locally definable continuous map between two locally definable manifolds with $e_X\in X$ and $e_Y\in Y$ such that $f(e_X)=e_Y$, then we have an induced homomorphism $f_*:\pi _1(X, e_X)\to  \pi _1(Y, e_Y):[\sigma ]\mapsto [f\circ \sigma ]$ with the usual functorial properties. \\
}
\end{defn}

We define the {\it o-minimal fundamental groupoid} $\Pi _1(X)$ of $X$ to be  the small category $\Pi _1(X)$  given by

\begin{displaymath}
\begin{array}{l}
\ob (\Pi _1(X))  = X, \\
\Ho _{\Pi _1(X)}(x_0, x_1)  = \raisebox{1ex}{$\mathsurround=0pt\displaystyle  {\mathbb P}(X, x_0, x_1)$}
\Big/ \raisebox{-1ex}{$\mathsurround=0pt\displaystyle \approx $}\\
\end{array}
\end{displaymath}
We  set $[\gamma ]:=$ the class of $\gamma \in {\mathbb P}(X, x_0,x_1)$. By Lemma~\ref{th-weak-homotopy-composition}, the small category $\Pi _1(X)$ is indeed a groupoid with operations 
$$\Ho_{\Pi _1(X)}(x_0,x_1)\times \Ho _{\Pi _1(X)}(x_1,x_2)\to \Ho _{\Pi _1(X)}(x_0,x_2)$$ 
given by $[\delta ] \circ [\gamma]=[\gamma \cdot \delta ]$.

Note that if $x\in X$, then ${\mathbb P}(X,x,x)={\mathbb L}(X,x)$ and so
$$\pi _1(X,x)=\Ho _{\Pi _1(X)}(X,x,x).$$

If $X$ is a locally definable manifold and $x\in X$, we define $\Pi _1(X,x)$ to be  the category given by
\begin{displaymath}
\begin{array}{l}
\ob (\Pi _1(X,x))  = \{x\}, \\
\Ho _{\Pi _1(X,x )}(x, x)  =  \pi _1(X,x). \\
\end{array}
\end{displaymath}

If $f:X\to Y$ is a locally definable continuous map between locally definable manifolds, then we have an induced functor $f_*:\Pi _1(X)\to \Pi _1(Y)$ which is a morphism of groupoids sending the object $x\in X$ to the object $f(x)\in Y$ and a morphism $[\gamma ]$ of $\Pi _1(X)$  to the morphism $[f\circ \gamma ]$ of $\Pi _1(Y)$. \\


\begin{lem}\label{lem Pi1 and x and connected}
Let $X$ and $Y$ be  locally definable manifolds. Then
\begin{enumerate}
\item
If $X$ is definably path connected then the natural functor $\Pi _1(X, x)\to \Pi _1(X)$ is an equivalence for every $x\in X$.
\item
The natural functor $\Pi _1(X\times Y) \to \Pi _1(X)\times \Pi _1(Y)$ given by projection  is an equivalence.
\end{enumerate}
\end{lem}

\pf
(1) The functor $\Pi _1(X, x)\to \Pi _1(X)$ sends the object $x$ of $\Pi _1(X,x)$ to the object $x$ of $\Pi _1(X)$ and sends a morphism of $\Pi _1(X,x)$ represented by a definable loop at $x$ to the morphim of $\Pi _1(X)$ represented by the same definable loop at $x$. By definition this morphism is fully  faithfull. Since $X$ is definably path connected,  every object of $\Pi _1(X)$ is isomorphic to the object $x$. So the functor is also essentially surjective. Therefore, it is an equivalence.

(2) 
The functor $\Pi _1(X\times Y) \to \Pi _1(X)\times \Pi _1(Y)$ sends a morphism of $\Pi _1(X\times Y)$ represented by a definable path $\rho $ in $X\times Y$ to the morphism of $\Pi _1(X)\times \Pi _1(Y)$ represented in each coordinate by  the definable paths $q_1\circ \rho $ in $X$ and $q_2\circ \rho $ in $Y$ where $q_1$ and $q_2$ are the projections onto $X$ and $Y$, respectively. This functor is an isomorphism with inverse given by the functor  $\Pi _1(X)\times \Pi _1(Y)\rightarrow   \Pi _1(X\times Y)$ that sends the object $(x,y)$ of  $\Pi _1(X)\times \Pi _1(Y)$ to the object $(x,y)$ of  $ \Pi _1(X\times Y)$ and sends  a  morphism  of $\Pi _1(X)\times \Pi _1(Y)$ represented by a pair of definable paths $\gamma $ in $X$ and $\delta $ in  $Y$  to the morphism  of  $ \Pi _1(X\times Y)$ represented by the definable path in $X\times Y$ with coordinates $\gamma $ and $\delta $. \qed \\

\begin{cor}\label{cor pi1 and x and connected}
Let $X$ and $Y$ be  locally definable manifolds with $e_X\in X$ and $e_Y\in Y$. Then
\begin{enumerate}
\item
If $X$ is definably path connected then $\pi _1(X, e_X)\simeq \pi _1(X, x)$ for every $x\in X$.
\item
$\pi _1(X, e_X)\times \pi _1(Y, e_Y)\simeq  \pi _1(X\times Y, (e_X,e_Y)).$
\end{enumerate}
\end{cor}



\medskip
\noindent
\textbf{Notation:}  As usual for a definably path connected locally definable manifold $X$ if there is no need to mention a base point $e_X\in X$, then by Corollary \ref{cor pi1 and x and connected} (1), we may denote $\pi _1(X,e_X)$ by $\pi _1(X)$.\\




\end{subsection}
\end{section}

\begin{section}{Topology on products of  definable group-intervals}\label{section top prod grp-int}
In this section we study some topology on products of  definable group-intervals including: definable normality, locally definable covering maps and the relativized new o-minimal fundamental group.

\begin{subsection}{Products of definable group-intervals}\label{subsection prod grp-int}
Here we recall a few notions about products of definable group-intervals. The results we will need came from \cite{epr} or are built  from what is done in that paper. \\

Recall the following (\cite[Definition 3.1]{epr}). Note that, as pointed out by Y. Peterzil to the authors,  in \cite[Definition 3.1 ]{epr}, condition  (ii) is only written down for {\it positive group-intervals} and it should be replaced by (ii) below in the case of group-intervals. Compare also with \cite[Defintion 2.8]{PL}.

\begin{defn}\label{defn gp-int}
{\em
A {\it definable group-interval}  $J=\langle (-b, b), 0, +,  <\rangle $ is an open interval $(-b, b)\subseteq M$, with $-b<b$ in $M\cup \{-\infty, +\infty \},$ together with a binary partial continuous definable operation $+:J^2\to J$ and an element $0\in J$, such that:
\begin{itemize}
\item[(i)]
 $x+y=y+x$ when defined; 
 $(x+y)+z=x+(y+z)$ when defined; 
 if  $x<y$ and $x+z$ and $y+z$ are defined then $x+z<y+z;$
\item[(ii)]
for every $x\in J,$ if $x>0,$ then the set $\{y\in J:\,\, x+y \,\, \textrm{is defined}\}$ is an interval of the form $(-b, r(x));$
\item[(iii)]
for every $x\in J,$ we have  $\lim _{z\to 0}(z+x)=x$ and  if $x>0$ we have also $\lim _{z\to r(x)^-}(x+z)=b;$
\item[(iv)]
for every $x\in J$ there exists $-x\in J$ such that $x+(-x)=0.$
\end{itemize}
The definable group-interval $J$ is {\it unbounded} (resp. {\it bounded}) if the operation $+$ in $J$ is total (resp. not total). 

The notion of a {\it definable homomorphism} between definable group-intervals is defined in the obvious way.
}
\end{defn}

\begin{nrmk}\label{nrmk unbounded}
{\em
Let  $J=\langle (-b, b), 0, +,  <\rangle $ be a definable group-interval. 
The following are equivalent: 
\begin{itemize}
\item
$J$ is  an unbounded definable group-interval; 
\item
 $J$ is  an ordered, abelian, divisible definable group (\cite[Theorem 2.1]{PiSte}); 
 \item
 for every $x\in J,$ if $0<x,$ then $r(x)=b.$
 \end{itemize}
}
\end{nrmk}

\begin{nrmk}
{\em
Let  $J=\langle (-b, b), 0, +,  <\rangle $ be a definable group-interval. Then: 
\begin{itemize}
\item
 the map $J\to J:x\mapsto -x$   is a continuous definable    bijection and we have $-0=0$, $-(-x)=x$ and $0<x$ if and only if $-x<0$; 
 \item
 for every $x\in J,$ if $x<0$, then the set $\{y\in J: \,\, x+y \,\, \textrm{is defined}\}$ is the interval  $(-r(-x), b);$ 
 \item
  for every $x\in J,$ if $x<0$, then   $\lim _{z\to -r(x)^+}(x+z)=-b.$ 
  \end{itemize}
}
\end{nrmk}

\begin{nrmk}\label{nrmk grp-int}
{\em
Let  $J=\langle (-b, b), 0, +,  <\rangle $ be a definable group-interval. Then:
\begin{itemize}
\item
For every $x\in J,$ if  $x>0,$ then $r(x)>0$ (by (i) and  (iii)). So  for every $x\in J,$ if $x>0$ we have $0+x=x=x+0.$ (Both sides are defined and they are equal). 
 
  \item
The map $r:(0,b)\to (0,b)$ is  definable and non-increasing, i.e., for all $x, y>0,$ if $x<y$ then $r(y)\leq r(x).$ 

 To see this consider the definable, continuous and strictly increasing (by (i)) map $T_z:(0, r(z))\to (-b,b): t\mapsto z+t.$  If we had $r(x)<r(y)$ then for all $t\in (0,r(x))$ we would have $T_x(t)<T_y(t)$ (by (i)) and so $\lim _{t\to r(x)^-}T_y(t)=b$ (by (iii)). But then for $r(x)<s<r(y),$ since $T_y(s)<b$ we could find $t<r(x)$ such that $T_y(t)>T_y(s)$ contradicting the fact that $T_y$ is strictly increasing. 

 \item
 For all $x,y>0$ if $x+y$ is defined then $y+x$ is defined. Suppose not. Then $y<r(x)$ but $r(y)\leq x.$ Since $\lim _{z\to r(y)}y+z=b,$ by  (iii), there is $t<r(y)$ such that $x+y<y+t.$ Since $t<x$ and $r$ is  non-increasing we have $r(x)\leq r(t).$ So $y<r(t)$ and $t+y$ is defined. Since both $y+t$ and $t+y$ are defined, they are equal by (i). Thus $x+y<t+y$ but, since $t<x,$ we also have by (i) $t+y<x+y$ which is a contradiction.
 
 \item
 For all $x,y,z>0$ if $(x+y)+z$ is defined then $x+(y+z)$ is defined.  Suppose not, in particular $(y+z)+x$ is not defined. Then $z<r(x+y)\leq r(x), r(y)$ ($r$ is non-increasing and $x,y<x+y$) and $r(y+z)\leq x<r(y)$ since $x+y$ is defined. Since $\lim _{t\to r(y+z)}(y+z)+t=b,$ by  (iii), there is $t<r(y+z)$ such that $(x+y)+z<(y+z)+t=t+(y+z).$ Since $t<x<r(y)$ we have that  $y+t$ is defined and so $t+y$ is defined. Since $t<x$ we have $t+y<x+y$ and as $r$ is  non-increasing we have $r(x+y)\leq r(t+y).$ So $z<r(t+y)$ and $(t+y)+z$ is defined. Since both $t+(y+z)$ and $(t+y)+z$ are defined, they are equal by (i). Thus $(x+y)+z<(t+y)+z$ but, since $t<x,$ we also have by (i) $t+y<x+y$ and $(t+y)+z<(x+y)+z$ which is a contradiction.

\item
The map $r:(0,b)\to (0,b)$ is continuous.

To see this fix $x.$ Let $y<r(x).$  Let $z<y.$ So $x+z$ is defined and $T_z(0)<r(x).$ By continuity of $T_z,$ there is $0<\epsilon <r(x), r(z), r(x+z)$ such that for all $s\in (-\epsilon , \epsilon)$ we have $T_z(s)<r(x).$  So $z+s<r(x)$ and since  $s<r(z),$ we have that $z+s$ is defined and  $x+(z+s)$ is defined. But then by the above $s+z$ is defined and $x+(s+z)$ is defined and therefore by the above again $(x+s)+z$ is defined.  Thus $z<r(x+s)$ for all $z<y.$ So $y\leq r(x+s)$ for some $0<s<r(x).$ Since $r$ is non-increasing this is enough.\\
\end{itemize}
}
\end{nrmk}

By the (omitted arguments in the) proof of \cite[Lemma 3.5]{epr} we have:

\begin{fact}\label{fact grp-int grp by 4}
Let $J=\langle (-b, b), 0, +, -,  <\rangle $ is a definable group-interval. Then $J$ is divisible and there exists an injective, strictly increasing, continuous definable homomorphism $\tau :J\to J$ given by $\tau (x)=\frac{x}{4}$ such that if  $x, y\in \tau (J)=(-\frac{b}{4}, \frac{b}{4}),$ then $x+y,$  $x-y$ and $\frac{x}{2}$ are defined in $J.$ (Note that when $J$ is an unbounded definable group-interval, then $-\frac{b}{4}=-b$ and $\frac{b}{4}=b.$) \\
\end{fact}

\pf
By o-minimality,  let $c=\lim _{x\to b}r(x).$ Then, since $r$ is continuous non-increasing, $c\leq r(x)$ for all $x\in (0,b).$ 

Suppose $c=b.$ Then $r(x)=b$ for all $x$ in $(0,b).$ So by Remark \ref{nrmk unbounded}, $J$ is an ordered, abelian divisible group. Letting $\frac{b}{n}=b$ for every $n,$ the result follows.

Suppose that $c<b.$  By (iii) $\lim_{x\to 0^+}r(x)=b.$ So since $r$ is continuous non-increasing, there is a unique point in $(0,b),$ which we shall denote $\frac{b}{2},$ such that $r(\frac{b}{2})=\frac{b}{2}$ and for all $x<\frac{b}{2}$ we have $x<r(x).$ In particular, for all $x<\frac{b}{2},$ $2x=x+x$ is defined. Also we have that the definable map $(0,\frac{b}{2})\to (0,b):x\mapsto 2x$ is continuous, strictly increasing, $\lim _{x\to 0^+}2x=0$ and  $\lim _{x\to \frac{b}{2}^-}2x=b$ (by (iii)). Therefore, $(-\frac{b}{2}, \frac{b}{2})$ with $0$, $+$ and $<$  restricted is a definable subgroup-interval of $J$ and   $(-\frac{b}{2}, \frac{b}{2})\to (-b,b):x\mapsto 2x$ is a continuous definable isomorphism. Moreover, if  $x, y\in (0, \frac{b}{2}),$ say $x<y,$ then $r(y)\leq r(x)$ and since $y<r(y)$ we have that $x+y$ is defined. Since for $x\in (0,\frac{b}{2})$ and $y\in (-\frac{b}{2},0]$ we have by (i) and (iii) that $x+y$ is defined, it follows that for all $x,y\in (-\frac{b}{2}, \frac{b}{2})$ both $x+y$ and  $x-y$  are defined in $J.$  

By induction, it follows that for every $n$ there is a point   in $(0,b),$ which we shall denote $\frac{b}{n},$ such that $\lim _{x\to \frac{b}{n}^+}nx=b,$  $(-\frac{b}{n}, \frac{b}{n})$ with $0$, $+$ and $<$  restricted is a definable subgroup-interval of $J$ and   $(-\frac{b}{n}, \frac{b}{n})\to (-b,b):x\mapsto nx$ is a continuous definable isomorphism (in particular, $J$ is divisible). Moreover,  for all $x,y\in (-\frac{b}{n}, \frac{b}{n})$ both $x+y$ and  $x-y$  are defined in $J.$  \qed\\

\begin{nrmk}\label{nrmk bJ struct}
{\em
Associated to $\bJ$ there is a definable  o-minimal structure ${\mathbb J}$ whose domain is obtained by ordering  disjoint copies of the $J_i$'s according to their numbering  and adding a point between two successive copies.  The ${\mathbb J}$-definable sets are those that are already definable in ${\mathbb M}.$ By \cite[Fact 4.5]{epr}, ${\mathbb J}$ has definable choice.\\
}
\end{nrmk}

Later we will need the following:\\

\begin{defn}
{\em
Let $\bJ=\Pi _{i=1}^mJ_i$  be a cartesian product of   definable group-intervals $J_i=\langle (-_ib_i, b_i), $ $0_i, +_i, -_i,  <\rangle $.  
\begin{itemize}
\item
We will say that a set $X$ is a {\it $\bJ$-set} if $X\subseteq \bJ;$ In particular, a {\it $\bJ$-cell} is a cell which is  a $\bJ$-set.
\item
We will that $X$ is a (locally) definable manifold with  definable {\it $\bJ$-charts} if $X$ has definable charts  $\{(U_l, \phi _l)\}_{l\leq \kappa }$ with each $\phi _l(U_l)$ a definable  $\bJ$-set.
\item
We will say that a set $X$ is a {\it $\bJ$-bounded set} if  $X\subseteq \Pi _{i=1}^m[-_ic_i,c_i]$ for some $c_i>0_i$ in $J_i;$ In particular, a {\it $\bJ$-bounded cell} is a cell which is  a $\bJ$-bounded set.
\item
We will that $X$ is a (locally) definable manifold with  definable {\it $\bJ$-bounded charts} if $X$ has definable charts  $\{(U_l, \phi _l)\}_{l\leq \kappa}$ with each $\phi _l(U_l)$ a definable  $\bJ$-bounded set.\\
\end{itemize}
}
\end{defn}

\end{subsection}

\begin{subsection}{Definable normality in  products of definable group-intervals}\label{subsection normal affine}
Here we study the notion of definably normal in products of definable group-intervals extending what was known in o-minimal expansions of ordered groups (\cite[Chapter 6, \S3]{vdd}).\\

Recall that a  definable space $X$  is \textit{definably normal} if one of the following equivalent conditions holds:
\begin{enumerate}
\item
for every disjoint closed definable subsets $Z_1$ and $Z_2$ of $X$ there are disjoint open definable subsets $U_1$ and $U_2$ of $X$ such that $Z_i\subseteq U_i$ for $i=1,2.$

\item
for every  $S\subseteq X$ closed definable  and $W\subseteq X$ open definable such that $S\subseteq W$, there is an open definable subsets $U$  of $X$ such that $S\subseteq U$ and $\bar{U}\subseteq W$.\\
\end{enumerate}

Definable normality is quite useful since it gives the shrinking lemma (compare with \cite[Chapter 6, (3.6)]{vdd}):\\

\begin{fact}[The shrinking lemma]\label{fact shrinking lemma}
Suppose that $X$ is a definably normal definable  space. If $\{U_i:i=1,\dots ,n\}$ is a covering of $X$ by open definable subsets, then there are definable open subsets $V_i$ and definable closed subsets $C_i$ of $X$ ($1\leq i\leq n$) with $V_i\subseteq C_i\subseteq U_i$ and $X=\cup \{V_i:i=1,\dots, n\}$. \\
\end{fact}

\begin{defn}\label{defn orth}
{\em
Two definable intervals $(b,b')\subseteq M$  and $(a, a')\subseteq M$ are {\it non-orthogonal} if   there  are sub-intervals  $(c, c')\subseteq  (b, b')$ and $(d, d')\subseteq  (a, a')$ together with a definable bijection $\sigma : (c, c')\to  (d, d').$ The intervals are {\it orthogonal} if they are not non-orthogonal.

Two be cartesian products   of definable open intervals, ${\mathbf J}=\Pi _{i=1}^mJ_i$ and ${\mathbf I}=\Pi _{j=1}^nI_j,$  are {\it orthogonal} if  for any $l\in \{1, \ldots , m\}$ and any $k\in \{1, \ldots , n\}$ we have that $J_l$ is orthogonal to $I_k.$\\
}
\end{defn}

This notion of orthogonality is a slight extension of the notion that appears in \cite{PePiSt00a}: here we do not assume that the open intervals are ``transitive''. There is a general notion of orthogonality of definable subsets of a (monster) model of a complete first order theory, see \cite[Section 2]{BerMam}. Our Lemma \ref{lem open and closed in orth} below shows that orthogonality as in Definition \ref{defn orth} coincides with the general notion of orthogonality.\\

\begin{lem}\label{lem map orth}
Let $\bI$ be a cartesian product of definable open intervals which is orthogonal to a definable open interval $J.$ Then every definable map $f:\bI \to J$ has finite image.
\end{lem}

\pf
We prove the result by induction on $n.$ For $n=1$ the result follows from orthogonality and the monotonicity theorem (\cite[Chapter 3, (1.2)]{vdd}).  Suppose that the result holds for $n-1.$ For each $t\in I_n$ let $f_t:\Pi _{j=1}^{n-1}I_j\to J$ be given by $f_t(x)=f(x,t).$ By the induction hypothesis, each $f_t$ has finite image. By uniform finiteness property \cite[Chapter 3, (2.13)]{vdd}, after replacing $I_n$ by a subinterval we may assume that for all $t\in I_n$ the image of $f_t$ has $p$ elements.   But then,  for each $k=1,\ldots , p$ we have a definable map $g_k:I_n\to J$ giving the $k$-th element of the image of $f_t.$ By orthogonality and the monotonicity theorem (\cite[Chapter 3, (1.2)]{vdd}) it follows that each $g_k$ has finite image, and so $f$ also has finite image. \qed \\

Below  and later, as usual, for $B\subseteq M^{m-1}$ and $f,g :B\to M$ we set
$$\Gamma (f)=\{(x,y)\in M^{m-1}\times M: x\in B\,\,\textrm{and}\,\,y=f(x)\},$$  
$$(f,g)_B=\{(x,y)\in M^{m-1}\times M: x\in B\,\,\textrm{and}\,\,f(x)<y<g(x)\}.$$  

\begin{lem}\label{lem open and closed in orth}
Let $\mathbf I$ and $\mathbf J$ be orthogonal cartesian products of definable open  intervals.  Let $A\subseteq \bI \times \bJ$ be a definable set and consider the  uniformly definable family  $\{A_x: x\in \bI\}$ of definable subsets $A_x=\{y\in \bJ: (x,y)\in A\}$ of $\bJ$. Then there are $x_0, \ldots , x_s\in \bI$ such that  $\{A_x: x\in \bI\}=\{A_{x_0},\ldots , A_{x_s}\}$. Moreover we have 
$$A=\bigcup _{i=0}^s(\{v\in \bI:A_v= A_{x_i}\}\times A_{x_i})$$
and  if $A$ is open, then $\{v\in \bI:A_v\supseteq A_{x_i}\}$ and  $A_{x_i}$ are open definable sets for each $i=0,\ldots , s.$

In particular, for every closed definable subset $B\subseteq \bI \times \bJ$ there are finitely many closed definable subsets $C_1, \ldots , C_p\subseteq \bI$ and $D_1,\ldots , D_q\subseteq \bJ$ such that 
$$B=\bigcup \{C_i\times D_j: i=1, \ldots , p \,\,\textrm{and}\,\,j=1,\ldots , q\}.$$
\end{lem}

\pf
 Take a cell decomposition ${\mathcal C}$ of $A.$ Then $\{C\cap A_x:C\in {\mathcal C}\}$ is an induced uniform cell decomposition on each fiber $A_x.$ We can clearly replace  $A$ by each $C\in {\mathcal C}.$  So we may suppose that $A$ is a cell and so  each $A_x$ is a cell of a fixed dimension.
 
 Let $E$ be the projection of $A$ into  $\bI.$  We define definable maps $a^i,b^i:\bI\to  J_i$ (possibly $a^i=b^i$) associated to $A$ and points $a^i(x),b^i(x)\in J_i$ (possibly $a^i(x)=b^i(x)$) associated to $A_x$  by induction on $m,$ in the following way.  If $m=1,$ then either $A=\Gamma (h)$ or $A=(f,g)_E$ for some continuous definable maps $h,f,g:E\to J_1$ with $f<_1g.$ In the first case we set $a^1=b^1=h$ and in the second case we set $a^1=f$ and $b^1=g.$ Suppose that for  every cell $C\in \Pi _{i=1}^{m-1}J_i$ we have defined definable maps $c^i,d^i:\bI\to  J_i$ for $i=1, \ldots, {m-1}$ associated to $C.$ Let $\pi (A)$ be the projection of $A$ into $\Pi _{i=1}^{m-1}J_i,$ then  either $A=\Gamma (h)$ or $A=(f,g)_{\pi (A)}$ for some  continuous definable maps  $h,f,g:\pi (A)\to J_m$ with $f<_mg.$  Let $p(x)=\langle x, a^1(x), \ldots , a^{m-1}(x)\rangle , q(x)=\langle x, b^1(x), \ldots , b^{m-1}(x)\rangle \in  \pi (A_x)\subseteq \bI \times \Pi _{i=1}^{m-1}J_i.$ We define $a^m,b^m:\bI\to J_m$ by setting $a^m(x)=h(p(x))$ and $b^m(x)=h(q(x))$ in the first case or $a^m(x)=f(p(x))$ and $b^m(x)=g(q(x))$ in the second case.
 
 We proceed with the proof of the Lemma  by induction on $m$. Suppose $m=1.$   If $A_x$ varies infinitely and definably with $x,$ then  least one of the definable maps $a^1,b^1: \bI\to  J_1$ associated to $A$ must have infinite image, contradicting orthogonality (Lemma \ref{lem map orth}). 
 
 For the case $m>1$, assume the result fails and let  $A$ as above  be such that  $A_x$ varies infinitely and definably as $x$ does. By induction, we know there are only finitely many sets of the form $\pi(A_x)$, where $\pi$ is projection to the first $m-1$ coordinates. Fix a $\pi(A_x)$ such that there are infinitely many distinct $A_x$ with this projection, and restrict to this family. Then the definable maps $a^i,b^i:\bI\to  J_i,$ $i=1,\ldots, m-1,$ associated to $\pi (A)$ are constant, but one of the definable maps $a^m,b^m:\bI\to  J_m$ associated to $A$ would have  infinite image, contradicting orthogonality (Lemma \ref{lem map orth}).
 
 
 Clearly we have $A=\bigcup _{i=0}^s(\{v\in \bI : A_v= A_{x_i}\}\times A_{x_i})$ and if $A$ is open  then each fiber  $A_{x_i}$ is also open. 
 
 Suppose that $A$ is open but some  $\{v\in \bI:A_v\supseteq  A_{x_i}\}$ is not open. Fix $v\in \{v\in \bI:A_v\supseteq  A_{x_i}\}$ such that for every open box $B$ around $v$ in $\bI$, $B\nsubseteq \{v\in \bI:A_v\supseteq  A_{x_i}\}.$ Thus, we can find points $z$ as close as we like to $v$ such that $A_{x_i}\nsubseteq A_z$ for any such $z$.  Consider the family of definable sets $\{A_{x_i}\setminus A_z : z\in\bI \,\, \textrm{and}\,\, A_{x_i}\nsubseteq A_z\}$. By the first part of the lemma, there are only finitely many sets in this family, so there is one that occurs for $z$ arbitrarily close to $v$, say $A_{x_i}\setminus A_{z_0} $. Fix any point $y\in A_{x_i}\setminus A_{z_0}.$ Then for any open box $B$ containing $v$, we can find $z\in B$ with $y\notin A_z$. But then any box in $\bI\times \bJ$ around the point $\langle v,y \rangle \in A$ must contain a point not in $A$, namely $\langle z,y \rangle $ for such a $z$, contradicting that $A$ is open.
 \qed \\

 \begin{lem}\label{lem norm in orth}
Let $\bI _1$ and $\bI_2$ be orthogonal cartesian products of definable group-intervals and set $\bI=\bI _1\times \bI _2.$  Let $A\subseteq \bI _1$ be a definably normal definable subset and let $B \subseteq \bI _2$ be a definably normal definable subset. Then $A\times B\subseteq \bI$ is a definably normal definable subset. 
\end{lem}

\pf
Let $S, T\subseteq A\times B$ be  closed, disjoint definable subsets. Then by Lemma \ref{lem open and closed in orth}, $S=\cup \{S_{1i} \times  S_{2i}: i=1,\ldots , s\}$ with each $S_{1i}\subseteq A$ a closed (in $A$) definable subset and each $S_{2i}\subseteq B$ a closed (in $B$) definable subset. Similarly,   $T=\cup \{T_{1j} \times  T_{2j}: j=1,\ldots , t\}.$

First suppose $s=1.$  Since  $T$ is disjoint from $S$ we have that  each $T_{1j} \times  T_{2j}$ has empty intersection with $S$ and therefore, either  $T_{1j}$ has empty intersection with $S_{11}$ or $T_{2j}$ has empty intersection with $S_{21}.$ Suppose the first case holds.  Since $A\subseteq \bI _1$ is  definably normal and  $B \subseteq \bI _2$ is definably normal, there exist $U_{11}\subseteq A$ an  definable subset  open  in $A$ and $V_{1j}\subseteq B$ an  definable subset open in $B$, containing $S_{11}$ and $T_{1j}$ respectively, with empty intersection. Let $U_{21}\subseteq A$ be an arbitrary open (in $A$) definable subset and let  $V_{2j}\subseteq B$ be an arbitrary open (in $B$) definable subset. Then  the products $U_j=U_{11}\times U_{21}\subseteq A\times B$ and $V_j=V_{1j}\times V_{2j}\subseteq A\times B$  are definable subsets open in $A\times B$, with empty intersection and containing  $S$ and $T_{1j}\times T_{2j}$ respectively. Now take $U=\cap \{U_j: j=1,\ldots , t\}$ and  $V=\cup\{V_j:j=1,\ldots , t\}. $ Then $U$ and $V$ are  definable subsets open in $A\times B$,  with empty intersection and containing  $S$ and $T$ respectively.

If $s>1$, by the above, one can take  definable subsets  $U_i$ and $V_i$  open in $A\times B$,  with empty intersection and containing  $S_{1i}\times S_{2i}$ and $T$ respectively.  Let $U=\cup \{U_i: i=1,\ldots , s\}$ and  $V=\cap\{V_i:i=1,\ldots , s\}. $ Then $U$ and $V$ are definable subsets open in $A\times B$,  with empty intersection and containing  $S$ and $T$ respectively.
\qed \\

\begin{lem}\label{lem norm Jm}
Let ${\mathbf J}=\Pi _{i=1}^mJ_i$  be a cartesian product of  non-orthogonal  definable group-intervals $J_i=\langle (-_ib_i, b_i), $ $0_i, +_i, -_i,  <\rangle $. Then every definable subset of $\Pi _{i=1}^mJ_i$ is  definably normal.
\end{lem}

\pf
By the definable homeomorphisms $\tau _i:(-_ib_i,b_i)\to (-_i\frac{b_i}{4}, \frac{b_i}{4})$ of  Fact \ref{fact grp-int grp by 4} we only have to show that every definable subset of $\Pi _{i=1}^m(-_i\frac{b_i}{4}, \frac{b_i}{4})$ is  definably normal.

By non-orthogonality and o-minimality, for each $l, k\in \{1, \ldots , m\},$ consider sub-intervals  $(c_l, d_l)\subseteq  (-_lb_l, b_l)$ with $-_lb_l<c_l<d_l<b_l$ and $(c_k, d_k)\subseteq  (-_kb_k, b_k)$ with $-_kb_k<c_k<d_k<b_k$ together with a definable bijection $\sigma _{lk}: (c_l, d_l)\to  (c_k, d_k).$  

By o-minimality  we may assume that $\sigma _{lk}$ is continuous, strictly monotonic bijection. By the strictly increasing definable homeomorphisms $\tau _i:(-_ib_i,b_i)\to (-_i\frac{b_i}{4}, \frac{b_i}{4})$ of  Fact \ref{fact grp-int grp by 4} we can further assume that  $(c_l, d_l)\subseteq  (-_l\frac{b_l}{4}, \frac{b_l}{4})$ with $-_l\frac{b_l}{4}<c_l<d_l<\frac{b_l}{4}$ (resp.  $(c_k, d_k)\subseteq  (-_k\frac{b_k}{4}, \frac{b_k}{4})$ with $-_k\frac{b_k}{4}<c_k<d_k<\frac{b_k}{4}$). 

By Fact \ref{fact grp-int grp by 4}, composing $\sigma _{lk}$ with the translation $x\mapsto x+_l\frac{d_l+_lc_l}{2}$ on the right and with the translation $x\mapsto x-_k\frac{d_k+_kc_k}{2}$  on the left, we may assume that $c_l=-_ld_l$ and $c_k=-_kd_k.$ By Fact \ref{fact grp-int grp by 4} again, we still have $I_l=(-_ld_l, d_l)\subseteq  (-_l\frac{b_l}{4}, \frac{b_l}{4})$  with $-_l\frac{b_l}{4}<-_ld_l<d_l<\frac{b_l}{4}$(resp.  $I_k=(-_kd_k, d_k)\subseteq  (-_k\frac{b_k}{4}, \frac{b_k}{4})$  with $-_k\frac{b_k}{4}<-_kd_k<d_k<\frac{b_k}{4}$).

For each $i,$ set  $\sigma _i=\sigma _{i1}:I_i\to I_1.$ Composing each $\sigma _i$ with $(-_1\frac{b_1}{4}, \frac{b_1}{4})\to (-_1\frac{b_1}{4}, \frac{b_1}{4}):x\mapsto -_1x$ if necessary, we may assume that each $\sigma _i$ is a continuous  strictly increasing bijection.  For each $i$ let $\alpha _i:I_i\to (-_1b_1,b_1):x\mapsto \sigma _i(x)-\sigma _i(0_i).$ Then each $\alpha _i$ is a continuous strictly increasing injective map with $\alpha _i(0_i)=0_1.$

Consider the map 
$$\delta _i:(-_i\frac{b_i}{4}, \frac{b_i}{4})\times (-_i\frac{b_i}{4}, \frac{b_i}{4}) \to [0_1, \frac{b_1}{4})$$
  given by  $\delta _i(x,y) = \alpha _i(\min (|x-_iy|_i, \alpha _i^{-1}(d_1))).$ 
Since each $\alpha _i$ is continuous, each $\delta _i$ is continuous. 

Now let 
$$\delta :\Pi _{i=1}^m(-_i\frac{b_i}{4}, \frac{b_i}{4}) \times \Pi _{i=1}^m(-_i\frac{b_i}{4}, \frac{b_i}{4})\to [0_1, \frac{b_1}{4})$$
 be given by $\delta (u,v)=\max\{\delta _i(u_i,v_i):i=1,\ldots ,m \}$ whenever $u=\langle u_1, \ldots , u_m \rangle $ and $v=\langle v_1, \ldots , v_m \rangle .$  Then $\delta $ is a continuous  definable function. Moreover, if $A\subseteq B\subseteq \Pi _{i=1}^m(-_i\frac{b_i}{4}, \frac{b_i}{4})$ are definable subsets with  $A$ closed in $B$ and nonempty,  then the definable function 
 $$\delta _{A,B}:B\to [0_1, \frac{b_1}{4})$$
  given by $\delta _{A,B}(u)=\inf \{\delta (u,v):v\in A\}$ is continuous with $A=\{u\in B: \delta _{A,B}(u)=0_1\}.$
 Therefore, if $C$ and $D$ are nonempty, disjoint  definable subsets of a definable subset  $B\subseteq \Pi _{i=1}^m(-_i\frac{b_i}{4}, \frac{b_i}{4})$ which are closed in $B,$ then the definable subsets 
 $U=\{v\in B:\delta _{C,B}(v)<\delta _{D,B}(v)\}$ and $W=\{v\in B:\delta _{D,B}(v)<\delta _{C,B}(v)\}$
 of $B$ are open in $B$, disjoint and such that $C\subseteq U$ and $D\subseteq W.$

\qed \\

We are ready to prove the main observation of this subsection:

\begin{prop}\label{prop basis of norm}
Let $\bJ$ be a cartesian product of definable groups intervals. Then every open definable subset of $\bJ$ is a finite union of open  definably normal definable subsets.
\end{prop}

\pf
Let $\mathbf J_1,\ldots,\mathbf J_k$ be cartesian  products of non-orthogonal definable group-intervals, with $\bJ _i$ and $\bJ _j$ orthogonal for $i\ne j$, and $\bJ=\Pi_{i\le k}\bJ _i.$ We prove the result by  induction on $k.$ 

If $k=1$ then every open definable subset of $\bJ$ is  definably normal by Lemma \ref{lem norm Jm}. On the other hand, the inductive step follows from Lemmas \ref{lem open and closed in orth} and \ref{lem norm in orth}. 
\qed \\

The following consequence of Proposition \ref{prop basis of norm} will be useful later:

\begin{cor}\label{cor def jman basis norm}
Let $\bJ$ be a cartesian product of definable groups intervals. Suppose that $X$ is a locally definable manifold with  definable $\bJ$-charts.  If $Z$ is a locally closed definable subset of $X,$ then every  open definable subset of $Z$ is a finite union of open definably normal definable subsets.
\end{cor}

\pf
Intersecting $Z$ with the definable charts of $X$ we may identify $Z$ with a definable subset of $\bJ$ and the result follows from Proposition \ref{prop basis of norm}.
\qed \\

We end this subsection with some observations. We say that a   definable space $X$  is \textit{completely definably normal} if one of the following equivalent conditions holds:
\begin{enumerate}
\item
every definable subset $Z$ of $X$ is a definably normal definable subspace.

\item
every open definable  subset $U$ of $X$ is a definably normal definable subspace.
\item
for every closed definable subsets $Z_1$ and $Z_2$ of $X$, if $Z_0=Z_1\cap Z_2$, then  there are open definable subsets $V_1$ and $V_2$ of $X$ such that:
\begin{itemize}
\item[(i)]
$Z_i\setminus V_i=Z_0$, $i=1,2.$
\item[(ii)]
$V_1\cap V _2=\emptyset.$

\item[(iii)]
$\bar{V_1}\cap \bar{V_2}\subseteq Z_0.$
\end{itemize}

\item
for every definable subsets $S_1$ and $S_2$ of $X$, if $S_1\cap \bar{S_2}=\bar{S_1}\cap S_2=\emptyset $, then  there are disjoint open definable subsets $U_1$ and $U_2$ of $X$ such that $S_i\subseteq U_i$ for $i=1,2.$\\
\end{enumerate}

In topology a Hausdorff compact space is normal and moreover completely normal.  In the paper \cite[Theorem 2.11]{emp} it was shown that if ${\mathbb M}$ has definable choice, then every Hausdorff definably compact definable space is definably normal. The following shows that definable choice functions is not enough to guarantee  complete definable normality:\\

\begin{expl}\label{expl not comp norm}
{\em
Let $I_1=\langle (-_1d_1, d_1), $ $0_1, +_1, -_1,  <\rangle $ and $I_2=\langle (-_2d_2, d_2), $ $0_2, +_2, $ $-_2,  $ $<\rangle $  be  two orthogonal definable group-intervals. Let $(b_1, b'_1)\subseteq I_1$ and $(b_2, b'_2)\subseteq I_2$ be definable sub-intervals bounded in $I_1$ and $I_2$ respectively and consider also a point 
 $\langle a,b \rangle \in (b_1, b'_1)\times (b_2, b'_2).$ Take $X=[b_1, b'_1]\times [b_2, b'_2]$ and $U=((b_1, b'_1)\times (b_2, b'_2)) \setminus \{ \langle a,b \rangle  \},$ and take  $A=\{ a \} \times (b_2,b'_2)$ and $B=(b_1, b'_1)\times \{ b \}.$ 

Let $\bI =I_1\times I_2.$ Then $X$ is definable in a structure ${\mathbb I}$  with definable choice (Remark \ref{nrmk bJ struct}), $X$ is definably compact and  definably normal (Lemma \ref{lem norm in orth}) but the open definable subset $U$ of $X$ is not definably normal. Indeed, $A$ and $B$ are closed disjoint in $U$ which, by the description of open definable subsets of $\bI$ (Lemma \ref{lem open and closed in orth}), cannot be separated by disjoint open definable subsets. So $X$ is not completely definably normal.\\
}
\end{expl}

\end{subsection}

\begin{subsection}{Covering maps in  products of definable group-intervals}\label{subsection coverings}
Fix a cartesian product  ${\mathbf J}=\Pi _{i=1}^mJ_i$  of definable group-intervals $J_i=\langle (-_ib_i, b_i), $ $0_i, +_i, -_i,  <\rangle $.  Here we study the locally definable covering maps between locally definable manifolds with definable $\bJ$-charts  extending the results proved in o-minimal expansions of ordered group in \cite{eep}. \\

\begin{nrmk}
{\em
In this subsection we will  relativize to $\bJ$  the o-minimal fundamental group of Subsection \ref{subsection new fund grp}. 

One could possibly by-pass this by working in the o-minimal structure ${\mathbb J}$  (see Remark \ref{nrmk bJ struct}). In fact one of the authors, Marcello Mamino, has some notes  showing that Lemma \ref{lem pi of Jcell} below can be proved, with the o-minimal fundamental group of Subsection \ref{subsection new fund grp}, in arbitrary o-minimal structures with definable choice. However, in this other approach the proof of the analogue of Lemma \ref{lem pi of Jcell} below is  much longer and complicated. In any case we remark that the proofs of the analogues 
of the remaining results of this Subsection (Corollary \ref{cor pi cover  def jman}, Lemma \ref{lem ldef lifting paths} and Remark \ref{nrmk covers main}) would be exactly the same as the ones we present here.

There is yet another possibility of avoiding the need to relativize to $\bJ$ the notions of Subsection \ref{subsection new fund grp}. Since later in the applications we will work only with (locally) definable manifolds with $\bJ$-bounded charts, we could fix $0<c_i<\frac{b_i}{4}$ and assume that our (locally) definable manifolds have definable charts $\{(U_l,\phi _l)\}_{l\leq \kappa }$ with $\phi _l(U_l)\subseteq \Pi _{i=1}^m[-_ic_i,c_i].$ After that we would work in the definable o-minimal structure constructed by gluing (disjoint copies of) the intervals $[-_ic_i,c_i]$'s according to their numbering. This new o-minimal structure also has definable choice. In this case the proofs  would be that same as the ones below. The disadvantage would be the restriction on the locally definable manifolds considered.

Finally note that, since for any one of the  three possible o-minimal fundamental groups we have analogues of the results of this subsection,  by Remark \ref{nrmk covers main} in each one of the three situations we could argue as in \cite[Subsection 3.4]{eep} to conclude that these o-minimal fundamental groups are isomorphic on locally definable manifolds with definable $\bJ$-charts. Details will be carried out in another paper.\\
}
\end{nrmk}

\begin{defn}
{\em
$\,\,$
\begin{itemize}
\item
A basic $d$-$\bJ$-interval is a basic $d$-interval $\maI=\langle [a,b], \langle 0_{\maI}, 1_{\maI}\rangle \rangle$ with $[a, b]\subseteq J_l$ for some $l\in \{1, \ldots ,m\};$ a $d$-$\bJ$-interval is a $d$-interval $\maI=\maI _1\wedge\dotsb\wedge \maI _n$ with each $\maI _i$ a basic  $d$-$\bJ$-interval. Note that the $\maI_i$'s can be in different $J_l$'s.
\item
If $X$ is a locally definable manifold with definable $\bJ$-charts, then a definable $\bJ$-path (resp.  constant definable $\bJ$-path, or definable $\bJ$-loop) is a definable path (resp. constant definable path or definable loop) $\alpha :\maI \to  X$ with $\maI$ a $d$-$\bJ$-interval;  $X$ is definably $\bJ$-path connected if for every $u,v$ in $X$ there is a definable $\bJ$-path $\alpha :\maI \to   X$ such that $\alpha _0=u$ and $\alpha _1=v.$ 
\item
If $X$ and $Y$ are locally definable manifolds with definable $\bJ$-charts, then two definable continuous maps $f,g:Y\to X$ are definable $\bJ$-homotopic, denoted $f\sim _{\bJ} g,$ if there is a definable homotopy $F(t,s):Y\times \maJ\to   X$  between $f$ and $g$ with $\maJ$ a $d$-$\bJ$-interval;  two definable $\bJ$-paths $\gamma :\maI \to   X$, $\delta :\maJ\to   X$, with $\gamma _0=\delta _0 $ and $\gamma _1=\delta _1$, are definably $\bJ$-homotopic, denoted $\gamma \approx _{\bJ}\delta ,$ if there  are $d$-$\bJ$-intervals $\maI'$ and~$\maJ '$
such that $\maJ '\wedge\maI=\maJ \wedge\maI'$, and there is a definable $\bJ$-homotopy 
$${\mathbf c}_{\maJ '}^{\gamma_0}\cdot \gamma
\sim _{\bJ}
\delta\cdot  {\mathbf c}_{\maI'}^{\delta_1}$$
fixing the end points. 
\end{itemize}
}
\end{defn}

The results proved in Subsection \ref{subsection new fund grp} for the relations $\sim $ and $\approx $ hold also for $\sim _{\bJ}$ and $\approx _{\bJ}$ respectively. 

\begin{defn}
{\em
Let $X$ be a locally definable manifolds with definable $\bJ$-charts,  $e_X\in X$ and $x_0, x_1\in X$. Let ${\mathbb P}^{\bJ}(X, x_0, x_1)$ denote the set of all  definable $\bJ$-paths in $X$ that start at $x_0$ and end at $x_1$ and let  ${\mathbb L}^{\bJ}(X, e_X)$ denotes the set of all definable $\bJ$-loops that start and end at a fixed  element $e_X$ of $X$ (i.e. ${\mathbb L}^{\bJ}(X,e_X)={\mathbb P}^{\bJ}(X,e_X,e_X)$). Then the restriction of  $\approx _{\bJ}$   to ${\mathbb P}^{\bJ}(X, x_0, x_1)\times {\mathbb P}^{\bJ}(X, x_0, x_1)$ is an equivalence relation on ${\mathbb P}^{\bJ}(X, x_0, x_1)$ and
$$\pi _1^{\bJ}(X, e_X):=
\raisebox{1ex}{$\mathsurround=0pt\displaystyle {\mathbb L}^{\bJ}(X,e _X)$}
\Big/ \raisebox{-1ex}{$\mathsurround=0pt\displaystyle \approx _{\bJ}$}
$$
 is a group, the {\it o-minimal $\bJ$-fundamental group} $\pi _1^{\bJ}(X, e_X)$ of $X,$ with group operation given by $[\gamma ][\delta ]=[\gamma \cdot \delta ]$ and identity the class a of constant $\bJ$-loop at $e_X.$ Moreover, if $f:X\to  Y$ is a locally definable continuous map between two  locally definable manifolds with definable $\bJ$-charts with $e_X\in X$ and $e_Y\in Y$ such that $f(e_X)=e_Y$, then we have an induced homomorphism $f_*:\pi _1^{\bJ}(X, e_X)\to  \pi _1^{\bJ}(Y, e_Y):[\sigma ]\mapsto [f\circ \sigma ]$ with the usual functorial properties. \\
}
\end{defn}

\medskip
\noindent
\textbf{Notation:}  As usual for a definably $\bJ$-path connected locally definable manifold $X$ with definable $\bJ$-charts if there is no need to mention a base point $e_X\in X$, then by Corollary \ref{cor pi1 and x and connected} (1), we may denote $\pi _1^{\bJ}(X,e_X)$ by $\pi _1^{\bJ}(X)$.\\


We start with the following in which the proof of (1) is similar to that of  \cite[Chapter 6, Proposition (3.2)]{vdd}. \\ 

\begin{lem}\label{lem pi of Jcell}
 Let $C\subseteq \Pi _{i=1}^m(-_i\frac{b_i}{4}, \frac{b_i}{4})$ be a $\bJ$-cell. Then:
\begin{enumerate}
\item
 $C$ is definably $\bJ$-path connected. In fact there is a uniformly definable family of definable $\bJ$-paths connecting a given fixed point in $C$ to any other point in $C.$
 \item
 $C$ is definably $\bJ$-simply connected, i.e $\pi _1^{\bJ}(C)=1.$
 \end{enumerate}
\end{lem}

\pf
(1) We prove the result by  induction on the  definition of cells (recall that cells are defined inductively). If $C\subseteq \Pi _{i=1}^m\{0_i\}$  then the result  is clear. 

Suppose that   $C=\Gamma (h)$ for some continuous definable map $h:B\to (-_m\frac{b_m}{4}, \frac{b_m}{4})$  with $B\subseteq \Pi _{i=1}^{m-1}(-_i\frac{b_i}{4}, \frac{b_i}{4})\times \{0_m\}$ a cell. Then the projection of $C$ onto $B$ is a definable homeomorphism and (1) follows by induction hypothesis. 

Suppose that $C=(f,g)_B$ for some continuous definable maps $f,g:B\to (-_m\frac{b_m}{4}, \frac{b_m}{4})$ with $f<g$ and $B\subseteq \Pi _{i=1}^{m-1}(-_i\frac{b_i}{4}, \frac{b_i}{4})\times \{0_m\}$ a cell.  Moreover assume also that $-_{m}\frac{b_{m}}{4}<f<g<\frac{b_{m}}{4},$ if not the argument is easier. 

Let $\langle x, u\rangle , \langle x', u'\rangle \in C$ with $x, x'\in B.$  By Fact \ref{fact grp-int grp by 4} we can do the operations $\frac{f(z)+_{m}g(z)}{2}$ in the  component $J_{m}$ and so  $\langle x, \frac{f(x)+_{m}g(x)}{2}\rangle , \langle x', \frac{f(x')+_{m}g(x')}{2}\rangle \in C.$ Let $\alpha :\maI \to C$ be the vertical  definable $\bJ$-path with $\alpha _{0}=\langle x, u\rangle $ and $\alpha _{1}=\langle x, \frac{f(x)+_{m}g(x)}{2}\rangle $ and let $\alpha ':\maI ' \to C$ be the vertical  definable $\bJ$-path with $\alpha '_{0}=\langle x', \frac{f(x')+_{m}g(x')}{2}\rangle $ and $\alpha '_{1}=\langle x', u'\rangle .$ By the induction hypothesis, let $\beta : \maJ \to B$ be a definable $\bJ$-path with $\beta _{0}=x$ and $\beta _{1}=x'.$ Then $\alpha \cdot \gamma \cdot \alpha ': \maI \wedge \maJ \wedge \maI' \to C$ where $\gamma (t)=\langle \beta (t), \frac{f(\beta (t))+_{m}g(\beta (t))}{2}\rangle $ is a definable $\bJ$-path in $C$ connecting $\langle x, u\rangle $ to $ \langle x', u'\rangle .$ Since the definable $\bJ$-paths $\alpha, \alpha '$ and $\beta $ can be chosen uniformly, the definable $\bJ$-path $\gamma $ can also be defined uniformly.

(2) We prove the result by  induction on the  definition of cells. The first two cases are easy as above. Suppose that $C=(f,g)_B$ for some continuous definable maps $f,g:B\to (-_m\frac{b_m}{4}, \frac{b_m}{4})$ with $f<g$ and $B\subseteq \Pi _{i=1}^{m-1}(-_i\frac{b_i}{4}, \frac{b_i}{4})\times \{0_m\}$ a cell.  Moreover assume also that $-_{m}\frac{b_{m}}{4}<f<g<\frac{b_{m}}{4},$ if not the argument is easier. 

Let $\alpha :\maI \to C$ be a definable $\bJ$-loop at $\langle u,v\rangle .$ We want to show that $\alpha \approx _{\bJ}{\mathbf c}_{\maJ}^{\langle u,v\rangle }.$  By Fact \ref{fact grp-int grp by 4} we can do the operations $\frac{f(x)+_{m}g(x)}{2}$ in the  component $J_{m}$ and so  $\langle x, \frac{f(x)+_{m}g(x)}{2}\rangle \in C$ for all $x\in B.$ By (1) and  Corollary \ref{cor pi1 and x and connected} (1), we have  $\pi _1^{\bJ}(C,\langle u,v\rangle )\simeq \pi _1^{\bJ}(C, \langle u, \frac{f(u)+_{m}g(u)}{2}\rangle ).$ Therefore, after replacing $\alpha $ with the concatenation of a vertical definable $\bJ$-path with $\alpha $ we may assume that $\langle u,v\rangle =\langle u, \frac{f(u)+_{m}g(u)}{2}\rangle .$

Let $\pi :C\to B: \langle x,z\rangle \mapsto \langle x, 0_m\rangle $ be the projection. Let $\beta =\pi \circ \alpha :\maI \to B$ be the definable $\bJ$-loop at $u$ obtained by projecting $\alpha .$  Let $\gamma : \maI\to C$ the the definable $\bJ$-loop at $\langle u, \frac{f(u)+_{m}g(u)}{2}\rangle $ given by $\gamma (t)=\langle \beta (t),  \frac{f(\beta (t))+_{m}g(\beta (t))}{2}\rangle.$ Since $\gamma $ is a definable $\bJ$-loop in the cell $\Gamma (\frac{f+_mg}{2}),$ by the induction hypothesis, it is enough to show that  $\alpha  \approx _{\bJ}\gamma .$ By Remark \ref{nrmk sim and approx} it is enough to show that $\alpha \sim _{\bJ}\gamma .$

Let $\tau :C\to J_m: \langle x,z\rangle \mapsto z $ be the projection. Let $\mu :\maI \to C$ be the definable $\bJ$-path given by 
$$\mu (t)= \langle \beta (t), \min\,(\tau \circ \alpha (t), \tau\circ \gamma (t))\rangle .$$
Let $a=\max \{\max \,(\tau\circ \alpha (t), \tau \circ \gamma (t)):t\in \maI\}$ and let $b=\min \{\min \,(\tau\circ \alpha (t), \tau \circ \gamma (t)):t\in \maI\}.$  If $a=b$ then $\alpha =\gamma ,$ so  we may assume that $b<a.$ Consider the basic $d$-$\bJ$-interval $\maK=\langle [b,a], \langle a,b\rangle \rangle .$ Let $F:\maI\times \maK \to C$ be the continuous definable map given by 
$$F(t,r)=\langle \beta (t), \max \,(\tau \circ \mu (t), \min \, (\tau \circ \alpha (t),r))\rangle .$$
Then $F_0=\alpha $ and $F_1=\mu .$ Therefore, $\alpha \sim _{\bJ}\mu .$ Similarly, $\gamma \sim _{\bJ}\mu .$







\qed \\

By Theorem \ref{thm open cells} and  Lemma \ref{lem pi of Jcell} we have the following. Compare with the corresponding results \cite[Lemma 2.9 and Proposition 3.1]{eep} in o-minimal expansions of ordered groups.\\

\begin{cor}\label{cor pi cover  def jman}
Let $X$ be a definable  manifold of dimension $n$ with definable $\bJ$-charts. Then the following hold:
\begin{enumerate}
\item
$X$ is definably connected  if and only if $X$ is definably $\bJ$-path connected. In fact, for any definably connected definable subset $D$ of $X$ there is a uniformly definable family of definable $\bJ$-paths in $D$ connecting a given fixed point in $D$ to any other point in $D$.
\item
$X$ has  an  admissible cover  $\{O_s\}_{s\in S}$ by  open definably  connected definable subsets such that:
\begin{itemize}
\item
$\{O_s\}_{s\in S}$ refines the definable charts of $X$;

\item
for each $s\in S$, $O_s$  is definably homeomorphic to a  $\bJ$-cell of dimension $n,$ in particular,  the o-minimal $\bJ$-fundamental group $\pi _1^{\bJ}(O_s)$ is trivial.\\
\end{itemize}
\end{enumerate}
\end{cor}

We will need one further crucial result, Lemma \ref{lem ldef lifting paths} below. But first we need to recall a few definitions. See for example \cite{eep}.\\

\begin{defn}
{\em
Given a  definably connected locally definable manifold $S$, a locally definable manifold $X$ and an admissible cover ${\mathcal U}=\{U_{\alpha }\}_{\alpha \in I}$ of $S$ by open definable subsets, we say that a continuous surjective locally definable map $p_X:X\rightarrow S$ is a {\it locally definable covering map trivial over ${\mathcal U}=\{U_{\alpha }\}_{\alpha \in I}$} if the following hold:
\begin{enumerate}
\item[$\bullet $]
$p_X^{-1} (U_{\alpha })=\bigsqcup _{i\leq \lambda }U_{\alpha }^i$ a disjoint union of open definable subsets of $X;$

\item[$\bullet  $]
each $p_{X|U_{\alpha }^i}:U_{\alpha } ^i\rightarrow U_{\alpha }$ is a definable homeomorphism.
\end{enumerate}
A {\it locally definable covering map} $p_X:X\rightarrow S$ is a locally  definable covering map trivial over some admissible cover ${\mathcal U}=\{U_{\alpha }\}_{\alpha \in I}$ of $S$ by open definable subsets.

We say that two locally definable covering maps $p_X:X\rightarrow S$ and $p_Y:Y\rightarrow S$ 
are {\it locally definably homeomorphic} if there is a locally definable homeomorphism $F:X\rightarrow Y$ such that:

\begin{enumerate}
\item[$\bullet $]
$p_X=p_Y\circ F.$
\end{enumerate}

A locally definable covering map $p_X:X\to  S$ is {\it trivial} if it is locally definably homeomorphic to a locally definable covering map $S\times M\to S :(s,m)\mapsto s$ for some  set $M.$\\
}
\end{defn}

Let $p_Y:Y\to  T$ be a locally definable covering map, $X$ be a  locally definable manifold  and let $f:X\to  T$ be a locally definable  map.  A {\it lifting of $f$} is a continuous map $\tilde{f}:X\to  Y$  such that $p_Y\circ \tilde{f}=f$. Note that a lifting of a continuous  locally definable  map need not be a 
locally definable  map. However, if $X$ is definably connected, then any two continuous locally definable liftings which coincide in a point must be equal \cite[Lemma 2.8]{eep}.\\

 For analogues of the following lemma compare with \cite[Section 2]{eo} in o-minimal expansions of fields or with \cite[Lemma 2.13]{eep} in o-minimal expansions of ordered groups. In all three cases the proofs are the same, they only use the fact that the domains of the corresponding definable paths and definable homotopies are definably normal. \\

\begin{lem}\label{lem ldef lifting paths}
Let $X$ and $S$ be  locally definable  manifolds with definable $\bJ$-charts. Suppose that $p_X:X\to  S$ is a locally definable covering map. Then
the following hold.
\begin{enumerate}
\item
Let $\gamma :\maI  \to S$ be a definable $\bJ$-path in $S.$ Let $x\in X$be such that $p_X(x)=\gamma _{0}.$ Then there exists a unique definable $\bJ$-path $\tilde{\gamma }:\maI \to X$ in $X$ lifting $\gamma $ such that $\tilde{\gamma } _{0}=x$.

\item
Suppose that $F: \maI \times \maJ \to S$ is a definable $\bJ$-homotopy  between the definable $\bJ$-paths $\gamma $ and $\sigma $ in $S$. Let $\tilde{\gamma }$ be a definable $\bJ$-path in $X$ lifting $\gamma $. Then there exists a definable $\bJ$-path $\tilde{\sigma }$ in $X$ lifting $\sigma $ and there exists is a unique definable lifting $\tilde{F}: \maI \times \maJ  \to X$ of $F$, which is a definable $\bJ$-homotopy between $\tilde{\gamma }$ and $\tilde{\sigma }.$
\end{enumerate}
\end{lem}

\pf 
Let  ${\mathcal U}=\{U_{\alpha }\}_{\alpha \in I}$ be an admissible cover of $S$ by open definable subsets over which $p_X:X\to  S$ is trivial. We may assume that  ${\mathcal U}=\{U_{\alpha }\}_{\alpha \in I}$ refines the  definable charts of $S$ witnessing the fact that $S$ is a locally definable manifold with definable $\bJ$-charts.  

(1) First we assume that $\maI$ is a basic $d$-$\bJ$-interval $\langle [a,b], \langle 0_{\maI}, 1_{\maI}\rangle \rangle .$ We may also assume that the definable total order  $<_{\maI}$ on the domain $[a, b]$ of $\maI$ is $<.$ If not, the argument is similar, one just has to construct the lifting from right to left instead of from left to right.

Let $L\subseteq I$ be a finite subset such that $\gamma ([a, b])\subseteq {\bigcup_{l\in L}}U_l$. Then $[a,b]\subseteq  \bigcup_{l\in L}\gamma^{-1}(U_l)$, with  the $\gamma^{-1}(U_l)$'s open in $[a,b]$. Since $[a,b]$ is a Hausdorff, definably compact  definable space, definable in the o-minimal structure ${\mathbb J}$ with definable choice (Remark \ref{nrmk bJ struct}) it follows that $[a,b]$ is  definably normal (\cite[Theorem 2.11]{emp}). So by the shrinking lemma (Fact \ref{fact shrinking lemma}),  for each $l\in L$ there is $W_l\subset [a,b]$, open in $[a,b]$ such that $W_l\subset\overline{W_l}\subset\gamma^{-1}(U_l)$ and $[a,b]\subseteq  \bigcup_{l\in L} W_l$. Therefore,  there are $a=s_0<s_1<\cdots <s_r=b$ such that for each $i=0,\dots,r-1$ we have $\gamma([s_i,s_{i+1}]) \subset U_{l(i)}$ (and $\gamma(s_{i+1})\in U_{l(i)}\cap U_{l(i+1)}$). 

Lift $\gamma_1=\gamma_{|[a,s_{1}]}$ to $\tilde{\gamma_1}=(p_{|U_{l(0)}^{i_0}})^{-1}\circ \gamma _{|[a,s_1]}$, with  $\tilde{\gamma_1}_0=x$, using the definable homeomorphism $p_{|U_{l(0)}^{i_0}}\colon U_{l(0)}^{i_0}\to U_{l(0)}$, where $U_{l(0)}^{i_0}$ is the definable connected component of $p^{-1}(U_{l(0)})$ in which $x$ lays. Repeat the process for each $\gamma_{i+1}=\gamma_{|[s_i,s_{i+1}]}$ with $\tilde{\gamma_i}(s_i)$ instead of $x$. Patch the liftings together to obtain $\tilde{\gamma }.$ 

Now if $\maI =\maI _1\wedge \ldots \wedge \maI _k$ with each $\maI_i$ a basic $d$-$\bJ$-interval apply the previous process to lift $\gamma_1=\gamma_{|\maI _1}$ to $\tilde{\gamma_1}$, with  $\tilde{\gamma_1}_0=x$ and repeat the process for each $\gamma_{i+1}=\gamma_{|\maI _{i+1}}$ with $\tilde{\gamma_i}(1_{\maI_i})$ instead of $x$. Patch the liftings together to obtain $\tilde{\gamma }.$

Uniqueness follows (in each step) from \cite[Lemma 2.8]{eep}.

(2) First assume  that $\maJ$ is a basic $d$-$\bJ$-interval $\langle [c,d], \langle 0_{\maJ}, 1_{\maJ}\rangle \rangle .$ We may also assume that the definable total order  $<_{\maJ}$ on the domain $[c, d]$ of $\maJ$ is $<.$ If not, the argument is similar, one just has to construct the lifting from top to bottom  instead of from bottom  to top.

To proceed we also assume that $\maI$ is a basic $d$-$\bJ$-interval $\langle [a,b], \langle 0_{\maI}, 1_{\maI}\rangle \rangle .$ We may furthermore  assume that the definable total order  $<_{\maI}$ on the domain $[a, b]$ of $\maI$ is $<.$ If not the argument is similar, one just has to construct the lifting from right to left instead of from left to right. 

Let $L\subseteq I$ be a finite subset such that $F([a,b]\times [c,d])\subseteq {\bigcup_{l\in L}}U_l$. Then $[a,b]\times [c,d] \subseteq \bigcup_{l\in L}F^{-1}(U_l)$, with  the $F^{-1}(U_l)$'s open in $[a, b]\times [c,d]$.  Since $[a,b]\times [c,d]$  is a Hausdorff, definably compact  definable space, definable in the o-minimal structure ${\mathbb J}$ with definable choice (Remark \ref{nrmk bJ struct}) it follows that $[a,b]\times [c,d]$  is  definably normal (\cite[Theorem 2.11]{emp}).  Hence by the shrinking lemma (Fact \ref{fact shrinking lemma}),   we have that  for each $l\in L$ there is $W_l\subset [a, b]\times [c,d]$, open in $[a, b]\times [a,d]$ such that $W_l\subset\overline{W_l}\subset F^{-1}(U_l)$ and $[a, b]\times[c, d] \subseteq \bigcup_{l\in L}W_l$. Now take a cell decomposition of  $[a,b]\times [c,d]$ compatible with the $W_l$'s. This cell decomposition is given by  a decomposition  $a=t_0<t_1<\cdots <t_r=b$ of $[a, b]$ together with definable continuous functions $f_{i,j}:[t_i, t_{i+1}] \to [c,d]$ for $i=0,\ldots , r-1$ and $j=0,\ldots , k_i$ such that: (i) $f_{i,0}<f_{i,1}< \ldots <f_{i, k_i}$ for $i=0,\ldots , r-1;$ (ii) $\Gamma (f_{i,0})=[t_i, t_{i+1}]\times \{c\}$ and $\Gamma (f_{i, k_i})=[t_i, t_{i+1}]\times \{d\}$ for $i=0,\ldots , r-1;$ (iii) the two-dimensional $\bJ$-cells are of form $C_{i,j=}(f_{i,j}, f_{i, j+1})_{(t_i, t_{i+1})}.$ For each two-dimensional $\bJ$-cell $C_{i,j}$ and each $l(i,j)$ such that $C_{i,j}\subset W_{l(i,j)}$, we have $F(\overline{C_{i,j}})\subset U_{l(i,j)}$ and  for any two-dimensional $\bJ$-cells $C_{i,j}$ and $C_{i',j'}$ in $[a, b]\times [c,d]$,   and for each $l(i,j), l(i',j'),$ such that $C_{i,j}\subset W_{l(i,j)}$ and $C_{i',j'}\subset W_{l(i',j')}$ we also have  $F(\overline{C_{i,j}}\cap\overline{C_{i',j'}})\subset U_{l(i,j)}\cap U_{l(i',j')}$. 

Lift $F_{0,1}=F_{|\bar{C_{0,1}}}$ to $\tilde{F_{0,1}}=(p_{|U_{l(0,1)}^{i_{0,1}}})^{-1}\circ F_{|\bar{C_{0,1}}}$,  using the definable homeomorphism $p_{|U_{l(0,1)}^{i_{0,1}}}\colon U_{l(0,1)}^{i_{0,1}}\to U_{l(0,1)}$, where $U_{l(0,1)}^{i_{0,1}}$ is the definable connected component of $p^{-1}(U_{l(0,1)})$ in which $\tilde{\gamma }([t_0, t_1])$ lays. Repeat the process for each $F_{0, j+1}=F_{|\bar{C_{0, j+1}}}$ with 
$\tilde{F_{0, j}}(\Gamma (f_{0,j}))$ instead of $\tilde{\gamma }([t_0, t_1])$. Patch the liftings together to obtain $\tilde{F_{0}}:[t_0,t_1]\times [c,d]\to X$ a definable lifting of $F_{|[t_0,t_1]\times [c,d]}$ which is a definable $\bJ$-homotopy between $\tilde{\gamma }_{|[t_0,t_1]}$ and $\tilde{\sigma }_{|[t_0,t_1]}.$ Repeat the above process again but now for each $i=1, \ldots , r-1$, starting in each case with $\tilde{\gamma }([t_i, t_{i+1}])$ and obtain the liftings $\tilde{F_{i}}:[t_i,t_{i+1}]\times [c,d]\to X$ a definable lifting of $F_{|[t_i,t_{i+1}]\times [c,d]}$ which is a definable $\bJ$-homotopy between $\tilde{\gamma }_{|[t_i,t_{i+1}]}$ and $\tilde{\sigma }_{|[t_i,t_{i+1}]}.$ These liftings  patch together to give a  definable lifting $\tilde{F}:[a,b]\times [c,d]\to X$  of $F$ which is a definable $\bJ$-homotopy between $\tilde{\gamma }$ and $\tilde{\sigma }.$

Now if $\maI =\maI _1\wedge \ldots \wedge \maI _k$ with each $\maI_i$ a basic $d$-$\bJ$-interval apply the previous process to lift $F_1=F_{|\maI _1\times [c,d]}$ to $\tilde{F_1}$, with  $\tilde{F_1}(\maI _1, c)=\tilde{\gamma }(\maI _1)$ and repeat the process for each $F_{i+1}=F_{|\maI _{i+1}\times [c,d]}$ with $\tilde{\gamma }({\maI_{i+1}})$ instead of $\tilde{\gamma }(\maI _1)$. Then patch these liftings together to obtain a  definable lifting $\tilde{F}:\maI \times \maJ \to X$ of $F$ which is a definable $\bJ$-homotopy between $\tilde{\gamma }$ and $\tilde{\sigma }.$

Now if $\maJ =\maJ _1\wedge \ldots \wedge \maJ _k$ with each $\maJ_j$ a basic $d$-$\bJ$-interval apply the previous process to lift $F_1=F_{|\maI \times \maJ_1}$ to $\tilde{F_1}$, with  $\tilde{F_1}(\maI , 0_{\maJ_1})=\tilde{\gamma }(\maI )$ and repeat the process for each $F_{j+1}=F_{|\maI \times \maJ _{j+1}}$ with $\tilde{F_j}(\maI, 1_{\maJ _j})$ instead of $\tilde{F_1}(\maI , 0_{\maJ_1})$. To finish patch these liftings together to obtain a definable  lifting $\tilde{F}:\maI \times [c,d]\to X$ of $F$ which is a definable $\bJ$-homotopy between $\tilde{\gamma }$ and $\tilde{\sigma }.$

As above, uniqueness follows from  \cite[Lemma 2.8]{eep}.
\qed \\

We end by observing that all the main results from \cite{eep} about locally definable coverings maps and o-minimal fundamental groups in o-minimal expansions of ordered groups also hold  for locally definable coverings maps between locally definable  manifolds with definable $\bJ$-charts and o-minimal $\bJ$-fundamental groups.

\begin{nrmk}\label{nrmk covers main}
{\em
Let ${\bf P}$ be  the full subcategory of locally definable spaces in ${\mathbb M}$ whose objects are the locally definable  manifolds with definable $\bJ$-charts. Then in the category ${\bf P}$ the following hold:

\begin{itemize}
\item[(P1)]
\begin{itemize}
\item[(a)]
every object of ${\bf P}$ which is definably connected is uniformly definably $\bJ$-path connected;
\item[(b)]
given a locally definable covering map $p_X:X\to S$ in ${\bf P}$ then: (i) every definable $\bJ$-path $\gamma $ in $S$ has a unique lifting $\tilde{\gamma }$ which is a definable $\bJ$-path in $X$ with a given base point; (ii) every definable $\bJ$-homotopy $F$ between definable $\bJ$-paths $\gamma $ and $\sigma $ in $S$  has a unique lifting $\tilde{F}$ which is a definable $\bJ$-homotopy between the definable $\bJ$-paths $\tilde{\gamma }$ and $\tilde{\sigma }$ in $X$.
\end{itemize}
\item[(P2)]
Every object of ${\bf P}$ has  admissible covers by definably $\bJ$-simply connected, open definable subsets refining any admissible cover by open definable subsets.
\end{itemize} 

It follows, as observed in the Concluding remarks of the paper \cite{eep}, that with (P1) and (P2) above one proves in exactly the same way all  the main results of the paper \cite{eep}. 

In fact,  besides (P1) and (P2) (and their consequences) everything else that is required is, on the one hand, results from \cite{ejp2}, 
which hold in arbitrary o-minimal structures (and for locally definable spaces as well), and  on the other hand,  \cite[Chapter 6, (3.6)]{vdd}, which is used to notice that the domains of the ``good'' definable paths are definably normal. In our case here the good definable paths are the definable $\bJ$-paths and their domains are Hausdorff, definably compact  definable spaces (Lemma \ref{lem d-int space}), definable in the o-minimal structure ${\mathbb J}$ with definable choice (Remark \ref{nrmk bJ struct}) and so they are definably normal by \cite[Theorem 2.11]{emp}.

The fact that (P1) and (P2) are the only requirements needed to develop this kind of theory is somewhat not surprising. Indeed in topology, where we have  good notions of paths and homotopies with the lifting of paths and homotopies property, all one  needs is existence of such nice open covers as in (P2). In the o-minimal context (here and in \cite{eep}),  the role that (P1) (b) and (P2) play is  similar to the role the analogue properties play in topology.  However, (P2) is often used in combination with the results from \cite{ejp2} mentioned above  to get local definability. Also (P1) (a) is required essentially only once and  to get local definability (see \cite[Proposition 2.18]{eep}), the other places where it is used, it is used to replace definably connected by definably path connected. \\
 }
 \end{nrmk}

Due to Remark \ref{nrmk covers main}, in the rest of the paper, when needed, we will freely  use the results of \cite{eep} in our context of  locally definable coverings maps between locally definable  manifolds with definable $\bJ$-charts and o-minimal $\bJ$-fundamental groups.\\


\end{subsection}
\end{section}

\begin{section}{Orientability and degree theory}\label{section orient degree}
In this section we collect a couple of notions and results about o-minimal sheaf cohomology and orientability that we will be using later. We develop degree theory,  compute the  cohomology with definably compact supports of $\bJ$-bounded cells and obtain results about the orientation sheaf of definable manifolds with definable $\bJ$-bounded charts.

\begin{subsection}{Preliminary remarks on o-minimal cohomology and orientability}\label{subsection prelim rmks}
For the readers convenience we collect here a couple of notions and results about o-minimal sheaf cohomology and orientability that we will be using later. The full details on these can be found in one of the papers   \cite{ejp1} (also \cite{bf}), \cite{ep1} and \cite{ep3}. \\

Given a definable space $X,$ the {\it o-minimal site} $X_{\df}$ on $X$ is the category $\op (X_{\df})$ whose objects are open definable subsets of $X$ and whose morphisms are inclusions and the admissible covers $\cov (U)$ of $U\in \op (X_{\df})$ are covers by open definable subsets with finite subcoverings. Given a commutative ring with unit $A,$ we will denote by $\mod(A_{X_{\df}})$ the category of sheaves of $A$-modules on $X$. 

Associated to a definable space $X$ we have its {\it o-minimal spectrum} $\tilde{X}$ which is a spectral topological space whose points are the  ultrafilters  of definable subsets of $X$ (also called in model theory, types concentrated on $X$) equipped with the topology generated by the open subsets of the form $\tilde{U}$, where $U\in \op (X_{\df})$ (i.e. is an open definable subset of $X$).

The tilde functor  determines a morphism of sites
$$\nu _X:\tilde{X}\to X_{\df}$$
given by the functor
$$\nu _X^t:\op (X_{\df})\to \op (\tilde{X}):U\mapsto \tilde{U}.$$

\begin{thm}[\cite{ejp1}]\label{thm main iso on sheaves}
The inverse image of  $\nu _X:\tilde{X}\to X_{\df}$ determines  an isomorphism of categories
$$\mod (A_{X_{\df }}) \to \mod (A_{\tilde{X}}):F\mapsto \tilde{F},$$
where $\mod (A_{\tilde{X}})$ is the category of $A$-sheaves on the topological space $\tilde{X}$.
\end{thm}

By Theorem \ref{thm main iso on sheaves} to develop sheaf theory in the category of  definable spaces is equivalent to developing sheaf theory in the category of o-minimal spectral spaces. For instance, suppose that $X$ is a definable space and that  $\Phi $ is a {\it family of definable supports} on $X$ (i.e. a collection of closed definable subsets of $X$ such that: (i)  $\Phi $ is closed under finite unions and (ii) every closed definable subset of a member of $\Phi $ is in $\Phi$). (See \cite[Definition 2.9]{ep1}).  Then $\tilde{\Phi}$, the collection of all closed subsets of tildes of members of $\Phi$, is a family of supports on $\tilde{X}$ and we define the {\it o-minimal cohomology of $X$ with support $\Phi$ and coefficients the sheaf $F$} in $\mod (A_{X_{\df}})$ to be 
$$H_{\Phi }^*(X;F):=H^*_{\tilde{\Phi }}(\tilde{X};\tilde{F}).$$

For example, in the rest of the paper, we denote by $A_X$ the constant sheaf with value  (the commutative ring with unit)  $A$  on  the o-minimal site $X_{\df}$ on $X;$   $H^*(X; A_X)$ denotes the o-minimal cohomology (with support  the family of all closed definable subsets of $X$) and given a closed definable subset $D\subseteq X,$ $H_D^*(X; A_X)$ denotes the o-minimal cohomology with support  the family of all closed definable subsets of $D.$ \\

In the paper \cite{ep1} we used this approach to develop the theory of $\Phi$-supported sheaves, where $\Phi$ is {\it a family of definably normal supports}, namely a family of definable supports such that: (1) each element of $\Phi$ is definably normal, (2) for each $S \in \Phi$ and each open definable neighborhood $U$ of $S$ there exists a closed definable neighborhood of $S$ in $U$ which is in $\Phi.$  

As explained in the section on orientation and duality in the paper \cite{ep3} this theory from \cite{ep1} is enough to develop the notion of orientability for definable spaces such that $c,$ the family of all definably compact subsets of $X,$ is a family of definably normal supports. More precisely: 

\begin{defn}\label{defn orient sheaf}
{\em
Let $X$ be a definable space of dimension $n$  such that $c$ is a family of definably normal supports. Suppose that  for every   open definable subset $U\subseteq X$ there exists an admissible  (finite) cover  $U_1,\dots , U_l$ of $U$   such that, for each $i$, we have:
\begin{equation*}
H_c^p(U_i;{\mathbb Z}_X)=
\begin{cases}
{\mathbb Z} \qquad \textmd{if} \qquad p=n\\
\\
\,\,\,\,\,\,\,\,\,\,\,\,\,\,\,\,\,\,\,\,\,\,\,\,\,\,\,\,\,\,\,
\\
0\qquad \textmd{if} \qquad p\neq n. \\
\end{cases}
\end{equation*}
Then we say that $X$  {\it has  orientation sheaf} and the {\it orientation sheaf of $X$} is  the locally constant sheaf $\Or_{X}$ on $X_{\df}$ with sections
$$\Gamma (U;\Or_{X})\simeq \Ho (H_c^n(U;{\mathbb Z}_X),{\mathbb Z}).$$

We say that $X$ is {\it orientable} if there exists an isomorphism $\ZZ _X\simeq {\mathcal Or}_X$  in $\mod ({\mathbb Z}_{X_{\df}}).$ Such isomorphisms are called {\it  orientations of $X$}.  If $X$ is orientable, then the associated  {\it orientation class} $\mu _X\in \Gamma (X; \Or _X)$ is the section image of the section $1_X\in \Gamma (X;\ZZ _X)$ by the orientation.

We also define the {\it stalk}  $\Or _{X,x}$ of the orientation sheaf $\Or _X$ of $X$ at $x\in X$ by 
$$\Or _{X,x}=\lind {x\in U}\Gamma (U; \Or _X)$$
where the limit is over all open definable subsets  $U\subseteq X$  of $X$ such that $x\in U.$ \\
}
\end{defn}

Recall the following from \cite[Subsection 4.1]{eep}. 

\begin{defn}
{\em
A {\it locally definable ${\mathbb Z}$-covering map $p:X\to S$ trivial over  ${\mathcal U}$} is a locally definable covering map  $p:X\to S$ trivial over  ${\mathcal U}$ with  a continuous locally definable right action $X\times {\mathbb Z}\to X:(x,a)\mapsto x^a$ of ${\mathbb Z}$ on $X$ such that for each $s\in S$ there is an induced right action $p^{-1}(s)\times {\mathbb Z}\to p^{-1}(s)$ making the fiber $p^{-1}(s)$ a ${\mathbb Z}$-torsor.

We say that two locally definable ${\mathbb Z}$-covering maps $p_X:X\rightarrow S$ and $p_Y:Y\rightarrow S$ 
are {\it locally definably homeomorphic} if there is a locally definable homeomorphism $F:X\rightarrow Y$ such that:

\begin{enumerate}
\item[$\bullet $]
$p_X=p_Y\circ F.$
\item[$\bullet$]
For every $x\in X$ and $a\in {\mathbb Z},$ we have $F(x^a)=F(x)^a.$
\end{enumerate}

A locally definable ${\mathbb Z}$-covering map $p_X:X\to  S$ is {\it trivial} if it is locally definably homeomorphic to a locally definable ${\mathbb Z}$-covering map $S\times M\to S :(s,m)\mapsto s$ for some  ${\mathbb Z}$-torsor $M.$\\

}
\end{defn}
  
By \cite[Proposition 4.3]{eep} (see also \cite[Example 4.2 (2)]{eep}) we have an  equivalence between the category of locally constant $\ZZ _X$-sheaves on $X_{\df}$ and the category of locally definable $\ZZ$-covering maps of $X.$ Note that \cite[Proposition 4.3]{eep} is proved in o-minimal expansions of ordered groups but it only uses \cite[Lemma 2.1 (1)]{ejp2} which holds in arbitrary o-minimal structures. 

By this equivalence of categories (and its proof) we have:

\begin{fact}\label{fact or and lief cov}
Let $X$ be a definable space of dimension $n$  such that $c$ is a family of definably normal supports. Suppose that $X$ has an orientation sheaf. Then there is a canonically associated locally definable $\ZZ$-covering map
$$w_{\oO}:W_{\oO}\to X$$
where 
$$W_{\oO}=\bigsqcup _{x\in X}{\Or}_{X,x}\,\,\,\text{and}\,\,\, w_{\oO}(s_x)=x,$$
 such that $X$ is orientable  if and only if the locally definable $\ZZ$-covering map $w_{\oO}:W_{\oO}\to X$ is trivial. 
Moreover, if $\{U_j\}_{i\in J}$ is an admissible cover of $X$ by open definable subsets such that for each $j\in J$ the restriction ${\Or}_{X|U_j}\simeq {\mathbb Z}_{X|U_j},$ then $w_{\oO}:W_{\oO}\to X$ is a locally definable $\ZZ$-covering map trivial over ${\mathcal U}=\{U_j\}_{j\in J}$  with, for each $j\in J$
\begin{itemize}
\item
$w^{-1}_{\oO}(U_j)=\bigsqcup _{m\in {\ZZ}}U_j^m;$
\item
each $w_{\oO|U_j}:U_j^m\to U_j$ is a definable homeomorphism 
\end{itemize}
where for $m\in {\ZZ}$, $U_j^m=\{s_x:\Gamma (U_j; {\Or}_X)\simeq \ZZ:s \mapsto m,\,\,\,x\in U_j\}.$ \\
\end{fact}

Later we shall also need the following consequence of the Alexander duality, both proved in \cite{ep3}:

 \begin{fact}\label{fact alexander dual}
 Let $X$ be a definable space of dimension $n$  such that $c$ is a family of definably  normal supports.  Suppose that $X$ is orientable. Let $Z$ be a definably compact definable subset of $X$ with $l$ definably connected components. Then there is an isomorphism 
 $$H_Z^n(X;{\mathbb Z}_X)\simeq \Ho(H^0(Z;{\mathbb Z}_X),{\mathbb Z})\simeq {\mathbb Z}^l$$
 induced by the given orientation.  In particular, by excision, if $U$ is an open definable subset of $X$ such that $Z\subseteq U$, then we have an isomorphism 
$$H^n_Z(U;\ZZ_X) \simeq \ZZ ^l$$
compatible with the  inclusions of open definable neighborhoods of $Z$ in $X$.\\
 \end{fact}

\end{subsection}

\begin{subsection}{Degree theory}\label{subsection degrees}
Here we introduce degree theory for continuous definable maps between definable manifolds  of positive dimension  $n$ on which the definably compact subsets are a family of definably normal supports and are orientable. \\

Since the functor $\Ho (\bullet , \ZZ )$ on the category  of abelian groups when restricted to the subcategory of torsion free abelian groups is exact, below we will denote it by  $(\bullet )^{\, \vee}.$ In particular,  we will use quite often the fact that if $f:A\to B$ is an isomorphism of torsion free abelian groups, then $f^{\,\vee }: B^{\, \vee }\to A^{\, \vee}$ is also an isomorphism of abelian groups.\\

\begin{defn}\label{defn fund class}
{\em
Let $X$ be a  definable manifold of positive dimension $n$ such that $c$ is a   family of definably normal supports and  is orientable.  Let $Z$ be a definably compact definable subset with $l$ definably connected components and $U$  an open definable subset of $X$ such that $Z\subseteq U.$
 
We call the element $\zeta _Z\in H^n_Z(U;\ZZ_X)^{\,\vee}$ corresponding to $(1, \ldots , 1)\in \ZZ ^l,$ via the isomorphism of Fact \ref{fact alexander dual},  the {\it fundamental class around $Z$}. If $X$ is definably connected and definably compact, then we call $\zeta _X$ the {\it fundamental class of $X$.}\\ 
}
\end{defn}

\begin{nrmk}\label{nrmk fund class}
{\em
Let $X$ be as above. Let $Z, Z_1\subseteq Z_2$ be definably compact, definable subset of $X.$ Then:
\begin{enumerate}
\item
If $Z=\emptyset ,$ then $\zeta _Z=0.$
\item
$\zeta _{Z_1}$ is the image of $\zeta _{Z_2}$ under the homomorphism 
$$H^n_{Z_2}(U;\ZZ_X)^{\, \vee}\to H^n_{Z_1}(U;\ZZ_X)^{\, \vee}$$ 
induced by inclusion.
\item 
If $Z$ is definably connected, then $H^n_Z(U;\ZZ_X)^{\, \vee}\simeq \ZZ$ and $\zeta _{Z}\in H^n_Z(U;\ZZ_X)^{\, \vee}$ is a generator. 
\item
If $X$ is definably connected and definably compact, then the fundamental class $\zeta _X\in H^n(X;\ZZ_X)^{\, \vee}$  of $X$ corresponds to the orientation class $\mu _X\in \Gamma (X; \Or _X).$\\
\end{enumerate}
}
\end{nrmk}

\begin{defn}\label{defn degree}
{\em
Let $X$ and $Y$ be   definable manifolds of positive dimension $n$ on which $c$ is a   family of definably normal supports and  are orientable.   Let  $f:Y\to X$ be  a definable continuous map. Let $Z$ be a definable connected, definably compact, nonempty definable subset of $X$ such that $f^{-1}(Z)$ is a definably compact definable subset of $Y$. 

We call {\it degree of $f$ over $Z$} the unique element ${\rm deg}_Zf\in \ZZ$ such that  the image of the fundamental class around $f^{-1}(Z)$  under the map
$$(f^*)^{\, \vee}: H^n_{f^{-1}(Z)}(Y;\ZZ _Y)^{\, \vee}\to H^n_Z(X;\ZZ_X)^{\, \vee}.$$
is given by 
$$(f^*)^{\, \vee} (\zeta _{f^{-1}(Z)})={\rm deg}_Zf\zeta _Z.$$

If $Y$ is definably compact and $X$ definably connected, then ${\rm deg}f:={\rm deg}_Xf$ is called the {\it degree of $f$}. Note that ${\rm deg}_Zf=0$ if $f^{-1}(Z)=\emptyset .$
}
\end{defn}

The next lemmas  establish some basic properties of the degree. Their proofs are classical but  we include them for completeness.\\

\begin{lem}\label{homeomdeg}
Let $X$ and $Y$ be   definable manifolds of positive dimension $n$ on which $c$ is a   family of definably normal supports and  are orientable. Suppose that  $X$ is  definably compact. Then the following holds.

\begin{enumerate}
\item
Let $V$ be a definable open subset of $X$.  Let  $f\colon V\to X$ be  the inclusion map. Let  $Z$  be a nonempty  definably connected definably compact definable subset of $V$. Then,  ${\rm deg}_Zf=1$.

\item
Let  $f\colon Y \to X$ be a definable homeomorphism onto an open  definable subset of $X$. Suppose  $Z$ is a nonempty definably connected definably compact definable subset of $X$ such that  $f^{-1}(Z)$ is definably compact. Then, ${\rm deg}_Zf=\pm 1$. 
\end{enumerate}
\end{lem}

\pf
(1) The dual $H^n_Z(X;\ZZ_X)^{\, \vee }\to H^n_Z(V;\ZZ_X)^{\, \vee }$ of the excision isomorphism   is the inverse of $(f^*)^{\, \vee}: H^n_Z(V;\ZZ_X)^{\, \vee }\to H^n_Z(X;\ZZ_X)^{\, \vee }$ and hence 
$(f^*)^{\, \vee }( \zeta_Z)= \zeta_Z$ (we have identified $\zeta_Z$ with its image through the
dual of the excision isomorphism).

(2) By $H^n_Z(X;\ZZ_X)^{\, \vee }\simeq \ZZ$ and the dual 
$H^n_{Z}(f(Y);\ZZ_X)^{\, \vee } \stackrel{}{\cong} H^n_Z(X;\ZZ_X)^{\, \vee }$ of the excision,   the composition 
$$H^n_{f^{-1}(Z)}(Y;\ZZ_Y)^{\, \vee }\stackrel{(f^*)^{\, \vee }}{\cong}H^n_{Z}(f(Y);\ZZ_X)^{\, \vee }
\stackrel{}{\cong}
H^n_Z(X;\ZZ_X)^{\, \vee }$$ must take $\zeta_{f^{-1}(Z)}$ to $\pm \zeta_Z$.
\qed \\

\begin{lem} 
Let $X$ and $Y$ be   definable manifolds of positive dimension $n$ on which $c$ is a   family of definably normal supports and  are orientable.  Suppose that $X$  is  definably compact. Let $f: Y\to X$ be a definable continuous map. Let  $Z\subset Z_1$ be 
definably compact nonempty subsets of $X$ such that $Z$ is definably connected  and $f^{-1}(Z)$ and $f^{-1}(Z_1)$ are definably compact. Then, 
$$(f^*)^{\, \vee }:H^n_{f^{-1}(Z_1)}(Y;\ZZ_Y)^{\, \vee }\to H^n_{Z_1}(X;\ZZ_X)^{\, \vee }$$ 
takes $\zeta_{f^{-1}(Z_1)}$ into
$({\rm deg}_Zf)\,\zeta_{Z_1}.$ Moreover, if  $Z_1$ is also definably connected then  we have 
${\rm deg}_Z\,f={\rm deg}_{Z_1}\,f$. 
\end{lem}

\pf
Consider the following commutative diagram
$$
\xymatrix{
H^n_{f^{-1}(Z)}(Y;\ZZ_Y)^{\, \vee } \ar[d]^{(i^*)^{\, \vee }}\ar[r]^{(f^*)^{\, \vee }} & H^n_Z(X;\ZZ_X)^{\, \vee }  \ar[d]^{(j^*)^{\, \vee }}\\
H^n_{f^{-1}(Z_1)}(Y;\ZZ_Y)^{\, \vee } \ar[r]_{(f^*)^{\, \vee }} & H^n_{Z_1}(X;\ZZ_X)^{\, \vee },
}
$$
where ${(i^*)^{\, \vee }}$ and ${(j^*)^{\, \vee }}$ are induced by the respective inclusion maps. Chasing 
$\zeta_{f^{-1}(Z)}$ through the diagram gives 
$\zeta_{f^{-1}(Z)}\mapsto {\rm deg}_Zf\,\zeta_Z\mapsto {\rm deg}_Zf\,\zeta_{Z_1}$,
respectively 
$\zeta_{f^{-1}(Z)}\mapsto \zeta_{f^{-1}(Z_1)}\mapsto (f^*)^{\, \vee }(\zeta_{f^{-1}(Z_1)}).$
\qed \\

We have the following   useful particular case:

\begin{cor}\label{pointdeg}
Let $X$ and $Y$ be  definably  compact, definable manifolds of positive dimension $n$ on which $c$ is a   family of definably normal supports and  are orientable. Suppose that   $X$ is  definably connected. Let $f: Y\to X$ be a definable continuous map. Then ${\rm deg}_xf:={\rm deg}_{\{x\}}f= {\rm deg} f$, for any $x\in X$.
\end{cor}

Finally we will also need:

\begin{lem}\label{sumdeg} 
Let $X$ and $Y$ be   definable manifolds of positive dimension $n$ on which $c$ is a   family of definably normal supports and  are orientable.   Suppose that  $X$ is  definably compact. Let $f: Y\to X$ be a definable continuous map. Let $Z\subset X$ be a definably compact set such that $f^{-1}(Z)$ is definably compact. Suppose $Y=\sqcup_{\lambda=1}^mY_{\lambda}$ is such that each $Y_{\lambda}$ is an open definable subset of $Y$ and $f^{-1}(Z)=\sqcup _{\lambda =1}^mf^{-1}(Z)\cap Y_{\lambda}.$ Then  ${\rm deg}_Zf=\sum_{\lambda=1}^m {\rm deg}_Zf^{\lambda},$
where $f^{\lambda}=f_{|Y_{\lambda}}\colon Y_{\lambda}\to X$.
\end{lem}

\pf
Firstly, note that the $Z'_{\lambda}=f^{-1}(Z)\cap Y_{\lambda}$ are clopen in  $f^{-1}(Z)$ and hence definably compact, so it makes sense to speak about the fundamental class around $Z'_{\lambda}$. Then observe that if $i^{\lambda}: (Y_{\lambda}, Y_{\lambda}\setminus  Z'_{\lambda})\to (Y,Y\setminus f^{-1}(Z))$  is the inclusion map, then the following diagram commutes
$$
\xymatrix{
\oplus_{\lambda=1}^m H^n_{Z'_{\lambda }}(Y_{\lambda};\ZZ_Y)^{\, \vee } \ar[dr]_{\sum_{\lambda=1}^m (f^{\lambda *})^{\, \vee }}  \ar[r]^{(i^{\lambda *})^{\, \vee }}&  H^n_{f^{-1}(Z)}(Y;\ZZ_Y)^{\, \vee } \ar[d]^{(f^*)^{\, \vee }}\\
  & H^n_Z(X;\ZZ_X)^{\, \vee }.
} 
$$
Now,   for each $p\in f^{-1}(Z)$, consider the maps 
$$\oplus_{\lambda=1}^m H^n_{Z'_{\lambda }}(Y_{\lambda};\ZZ_Y)\stackrel{(i^{\lambda *})}{\to} H^n_{f^{-1}(Z)}(Y;\ZZ_Y)\stackrel{(i^{p*})}{\to} H^n_{\{p\}}(Y;\ZZ_Y)$$ 
and observe that
$((i^{p*})\circ (i^{\lambda *}))(\zeta_{Z'_{\lambda}})= \zeta _p$ (all the components of $\zeta_{Z'_{\lambda}}$ go to zero except the component $\zeta_{Z'_{\mu}}$  containing $p$ which  goes to $\zeta _p$). Hence, by unicity of the fundamental class
$(i^{\lambda *})(\zeta_{Z'_{\lambda}}) =\zeta_{f^{-1}(Z)}$. Therefore, by the above diagram, we have 
\begin{align*}
({\rm deg}_Zf)\zeta_Z & = (f^*)(\zeta_{f^{-1}(Z)})\\
& =(f^*)((i^{\lambda *})(\zeta_{Z'_{\lambda}}))\\
& =\sum_{\lambda=1}^m (f^{\lambda *} (\zeta_{Z'_{\lambda}}) \\
& = (\sum_{\lambda=1}^m{\rm deg}_Zf^{\lambda})\zeta_Z.
\end{align*} 
\qed \\

\end{subsection}

\begin{subsection}{Cohomology with definably compact supports of $\bJ$-bounded cells}\label{subsection coho jmcells}
Here we compute the o-minimal cohomology with definably compact supports of $\bJ$-bounded cells and obtain results about the orientation sheaf of definable manifolds with definable $\bJ$-bounded charts. \\

Below we let  ${\mathbf J}=\Pi _{i=1}^mJ_i$  be a cartesian product of  definable group-intervals $J_i=\langle (-_ib_i, b_i), $ $0_i, +_i, -_i,  <\rangle $.  \\



In the paper \cite[Section 4]{ep2} a  generalization to the cartesian product ${\mathbf J}=\Pi _{i=1}^mJ_i$ of definable group-intervals  of a construction done  in o-minimal expansions of ordered groups in \cite[Section 7]{bf} is presented. The construction is similar but, while in \cite{bf} one needs to consider only one parameter  in \cite{ep2} one considers $m$ parameters, one  for each $J_i.$  For the readers convenience we recall this construction in Theorem \ref{thm Uct} below as well as a couple of results that we will be using.  \\

Below, as usual, for $B\subseteq M^{m-1}$ and $f,g :B\to M$ we set
$$\Gamma (f)=\{(x,y)\in M^{m-1}\times M: x\in B\,\,\textrm{and}\,\,y=f(x)\},$$  
$$(f,g)_B=\{(x,y)\in M^{m-1}\times M: x\in B\,\,\textrm{and}\,\,f(x)<y<g(x)\}$$  
and 
$$[f,g]_B=\{(x,y)\in M^{m-1}\times M: x\in B\,\,\textrm{and}\,\,f(x)\leq y\leq f(x)\}.$$  \\

For (1), (2) and (3) of the next result compare with \cite[Lemma 4.11 (1) and (2) and Remark 4.14]{ep2}  and \cite[Lemma 7.1]{bf}.

\begin{thm}\label{thm Uct}
Let $C$ be a $\bJ$-bounded cell. Assume that $C\subseteq \Pi _{i=1}^m[-_ic_i, c_i]$ for some $0_i<c_i<\frac{b_i}{4}$ in $J_i.$ Then there exists a  definable family 
$$\{C_{t_1,\ldots , t_m} : 0_i<t _i<\frac{b_i}{4}, \,\, i=1, \ldots , m\}$$ 
of definable subsets of $C$ such that the following hold:
\begin{enumerate}
\item 
$C = \bigcup_{t_1,\ldots, t_m}C_{t_1,\ldots , t_m}$ where the union is over all $m$-tuples $t_1,\ldots , t_m$.
\item 
If  $0_i < t'_i < t_i$ for all $i=1, \ldots , m$, then  $C_{t_1,\ldots , t_m}\subset C_{t'_1,\ldots , t'_m}.$ 
\item
There is a point $p_C\in C$ such that for all $t_1,\ldots , t_m $, if   $c_i<t_i$ for all $i=1,\ldots , m,$  then  $C_{t_1,\ldots , t_m}=\{p_C\}.$ 
\item 
The $C_{t_1,\ldots , t_m}$'s are  closed and $\bJ$-bounded subsets of $C$ (in particular, by \cite[Theorem 2.1]{ps}, they are definably compact).
\item 
If $K\subseteq C$ is a definably compact definable subset, then there are  $t_1,\ldots , t_m$  such that $K\subseteq C_{t_1, \ldots , t_m}.$\\
\end{enumerate}
\end{thm}

\pf
The definable family $\{C_{t_1,\ldots , t_m} : 0_i<t _i<\frac{b_i}{4}, \,\, i=1, \ldots , m\}$ of definable  subsets  of $C$ is constructed by induction on $m.$ Recall that by Fact \ref{fact grp-int grp by 4} the group-interval operations $x-_iy$, $x-_iy$ and $\frac{x}{2}$ are all defined in each coordinate of $ \Pi _{i=1}^m(-_i\frac{b_i}{4}, \frac{b_i}{4}).$
\begin{itemize}
\item[(i)]
If $C$ is a singleton in $J_1$, we define $C_{t_1} = C.$ In this case, $C_{t_1}$ is clearly a closed and $\bJ$-bounded subset of $C.$
\item[(ii)]
If $C = (d, e)\subseteq J_1,$ then $C_{t_1} =[d +_1 \gamma ^1_{t_1}, e -_1\gamma ^1_{t_1}]$ where $\gamma ^1_ {t_1} = \min\{\frac{e-_1d}{2}, t_1\},$ (in this way $C_{t_1}$ is non empty). In this case, $C_{t_1}$ is clearly a closed and $\bJ$-bounded subset of $C.$
\item[(iii)]
Suppose that   $C=\Gamma (h)$ for some continuous definable map $h:B\to (-_m\frac{b_m}{4}, \frac{b_m}{4})$  with $B\subseteq  \Pi _{i=1}^{m-1}[-_ic_i, c_i]\times \{0_m\}$ a cell. Then the projection of $C$ onto $B$ is a definable homeomorphism.  By induction $B_{t_1, \ldots, t_{m-1}}$ is defined. We put $C_{t_1, \ldots, t_{m-1}, t_{m}}=\Gamma (h_{|B_{t_1,\ldots , t_{m-1}}}).$ By induction, $C_{t_1, \ldots, t_{m-1}, t_{m}}$ is  a closed and $\bJ$-bounded subset of $C.$
\item[(iv)]
Suppose that $C=(f,g)_B$ for some continuous definable maps $f,g:B\to (-_m\frac{b_m}{4}, \frac{b_m}{4})$ with $f<g$ and $B\subseteq \Pi _{i=1}^{m-1}[-_ic_i, c_i]\times \{0_m\}$ a cell.   By induction $B_{t_1,\ldots, t_{m-1}}$ is defined. We put $C_{t_1,\dots , t_{m-1}, t_{m}} = [f +_{m}\gamma ^{m}_{t_{m}}, g -_{m} \gamma ^{m}_{t_{m}}]_{B_{t_1,\ldots , t_{m-1}}},$ where $\gamma ^{m}_{t_{m}} := \min(\frac{g-_{m} f}{2}, t_{m}).$   By induction   $B_{t_1, \ldots, t_{m-1}}$ is  a closed and $\bJ$-bounded subset of $B.$ Also,  for each $x\in B_{t_1, \ldots, t_{m-1}}$, the fiber $(\pi _{|C})^{-1}(x)\cap C_{t_1,\dots , t_{m-1}, t_{m}}$ is closed and $\bJ$-bounded. Let $(x,y)\in C$ be an element in the closure of $C_{t_1,\dots , t_{m-1}, t_{m}}.$ Then $x\in B_{t_1,\ldots , t_{m-1}}$ and $(x,y)\in (\pi _{|C})^{-1}(x)\cap C_{t_1,\dots , t_{m-1}, t_{m}}\subseteq  C_{t_1,\dots , t_{m-1}, t_{m}}.$ So  $C_{t_1,\dots , t_{m-1}, t_{m}}$ is a closed and $\bJ$-bounded subset of $C.$ 
 \end{itemize}
 So (1), (2), (3) and (4) follow easily by construction. It remains to show (5). Note that if  $0_i < t'_i < t_i$ for all $i=1, \ldots , m$, then  by construction $C_{t'_1,\ldots , t'_m}$ contains the interior relative to $C$ of  $C_{t_1,\ldots , t_m}.$ Thus $C$ a directed union $\bigcup_{t_1,\ldots, t_m}U_{t_1,\ldots , t_m}$ of a definable family of  relatively open definable subsets. In particular, if $K\subseteq C$ is a definably compact subset, then $\{K\cap U_{t_1,\ldots , t_m} : 0_i<t _i<\frac{b_i}{4}, \,\, i=1, \ldots , m\}$ is a definable family of open definable subsets of $K$ (which is closed and bounded by \cite[Theorem 2.1]{ps}) with the property that every finite subset of $K$ is contained in one of the $K\cap U_{t_1,\ldots , t_m}.$ Therefore, by \cite[Corollary 2.2 (ii)]{PePi07},  there are  $t_1,\ldots , t_m$  such that $K\subseteq U_{t_1, \ldots , t_m}\subseteq C_{t_1, \ldots , t_m}.$
\qed \\

For the following compare also with \cite[Corollary 3.3]{bf}, in fact  the proof is the same one only has to use instead the  group-interval operations $x+_iy$, $x-_iy$ and $\frac{x}{2}$ in each coordinate of $ \Pi _{i=1}^m(-_i\frac{b_i}{4}, \frac{b_i}{4})$ (by Fact \ref{fact grp-int grp by 4}).

\begin{fact}\cite[Lemma 4.8]{ep2}\label{fact bf acyclic cells}
Let $C$ be a $\bJ$-bounded cell such  that $C\subseteq \Pi _{i=1}^m[-_ic_i, c_i]$ for some $0_i<c_i<\frac{b_i}{4}$ in $J_i.$ Then $C$ is  acyclic, i.e. $H^p(C;{\mathbb Z}_C)=0$ for $p>0$ and $H^0(C; {\mathbb Z}_C)={\mathbb Z}$.
\end{fact}

We also have the analogue of \cite[Lemma 7.1]{bf}:

\begin{fact}\cite[Lemma 4.11]{ep2}\label{fact bf cells}
Let $C$ be a $\bJ$-bounded cell of dimension $r$ such  that $C\subseteq \Pi _{i=1}^m[-_ic_i, c_i]$ for some $0_i<c_i<\frac{b_i}{4}$ in $J_i.$ Then 
\begin{enumerate}
\item If $0_i < t'_i < t_i$ for all $i=1, \ldots , m$, then the inclusion $C_{t_1,\ldots , t_m}
\subset C_{t'_1,\ldots , t'_m}$ induces an isomorphism 
$$H^{p}(C \backslash C_{t_1,\ldots , t_m}; {\mathbb Z}_C) \simeq  H^{p}(C \backslash C_{t'_1,\ldots , t'_m}; {\mathbb Z}_C).$$
\item
The o-minimal cohomology of $C \backslash C_{t_1,\ldots , t_m}$ is given by
\begin{equation*}
H^{p}(C \backslash C_{t_1,\ldots , t_m}; {\mathbb Z}_C) =
\begin{cases}
{\mathbb Z} ^{1+\chi _{1}(r)}\qquad \textmd{if} \qquad p\in \{0, r-1\}\\
\\
\,\,\,\,\,\,\,\,\,\,\,\,\,\,\,\,\,\,\,\,\,\,\,\,\,\,\,\,\,\,\,
\\
0\qquad \,\,\,\,\,\,\, \,\,\,\,\,\,\,\,\,\,\, \, \textmd{if} \qquad p\notin \{0, r-1\}
\end{cases}
\end{equation*}
where $\chi _1:{\mathbb Z}\to \{0,1\}$ is the characteristic function of the subset $\{1\}.$\\
\end{enumerate}
\end{fact}

\begin{nrmk}\label{nrmk compact closed}
{\em
Let $X$ be a definable manifold with definable $\bJ$-charts. Then $X$ is a Hausdorff  ${\mathbb J}$-definable manifold and ${\mathbb J}$ has definable choice. Therefore, by \cite[Corollary 2.8]{emp}, every definably compact definable subset of $X$ is a closed definable subset. Hence $c,$ the family of definably compact subsets of $X,$ is a family of definable supports and the o-minimal cohomology $H_c^*(X;{\mathbb Z}_X)$ with definably compact supports of $X$ is well defined.\\
}
\end{nrmk}


\begin{prop}\label{prop c-cohom of j-cell}
Let $C$ be a $\bJ$-bounded cell of positive dimension $r$ such  that $C\subseteq \Pi _{i=1}^m[-_ic_i, c_i]$ for some $0_i<c_i<\frac{b_i}{4}$ in $J_i.$ Then 

\begin{equation*}
H_c^{l}(C ; {\mathbb Z}_C) =
\begin{cases} 
{\mathbb Z} \qquad \textmd{if} \qquad l=r\\
\\
\,\,\,\,\,\,\,\,\,\,\,\,\,\,\,\,\,\,\,\,\,\,\,\,\,\,\,\,\,\,\,
\\
0\qquad \textmd{if} \qquad l\neq r.
\end{cases}
\end{equation*}
Moreover, the inclusion induces an isomorphism $H_{C_{t_1, \ldots , t_m}}^*(C;{\mathbb Z}_C)\simeq H_c^*(C;{\mathbb Z}_C)$ for every $t_1, \ldots , t_m.$

In particular, for every $u\in C$ the inclusion induces an isomorphism 
$$H_{\{u\}}^*(C;{\mathbb Z}_C)\simeq H_c^*(C;{\mathbb Z}_C).$$
\end{prop}



\pf For short, write $ C_{\bar{t}}$ instead of $C_{t_1,\ldots , t_m}.$ We have the long exact sequence
$$
\xymatrix{
 H^i_{C_{\bar{t}}}(C) \ar[r]  &  H^{i}(C )  \ar[r] &H^i(C\setminus C_{\bar{t}} ) \ar[r]  & H^{i+1}_{C_{\bar{t}}}(C)  \ar[r]  &\\
}
$$
(where we omitted the coefficients) given by the  pair $(C, C\setminus C_{\bar{t}}).$
 
On the other hand,  the sections $\Gamma _{C_{\bar{t}}}(C, {\mathbb Z}_C)$ of ${\mathbb Z}_C$ on $C$ with support on $C_{\bar{t}}$ are by definition the kernel of the restriction morphism
$$\Gamma (C, {\mathbb Z}_C)\to \Gamma (C\setminus C_{\bar{t}}, {\mathbb Z}_C)$$
which is clearly injective since $C$ is an open neighborhood of $C_{\bar{t}}.$ Therefore, 
\begin{eqnarray*}
H^{0}_{C_{\bar{t}}}(C; {\mathbb Z}_C)&\simeq &\Gamma _{C_{\bar{t}}}(C,{\mathbb Z}_C)=0.\\
\end{eqnarray*}
So by this and Facts \ref{fact bf acyclic cells} and \ref{fact bf cells} (2), we have:

Case: $r=1.$

$$
\xymatrix{
0 \ar[r]  &  {\mathbb Z}  \ar[r] & {\mathbb Z}^2  \ar[r]  & H^{1}_{C_{\bar{t}}}(C; {\mathbb Z}_C)  \ar[r]  & 0\\
}
$$
and so, by exactness, in this case we  have 
\begin{equation*}
H_{C_{\bar{t}}}^{l}(C ; {\mathbb Z}_C) =
\begin{cases} 
{\mathbb Z} \qquad \textmd{if} \qquad l=r\\
\\
\,\,\,\,\,\,\,\,\,\,\,\,\,\,\,\,\,\,\,\,\,\,\,\,\,\,\,\,\,\,\,
\\
0\qquad \textmd{if} \qquad l\neq r.
\end{cases}
\end{equation*}

Case: $r>1.$

We have 
$$
\xymatrix{
0 \ar[r]  &  {\mathbb Z}  \ar[r] & {\mathbb Z}  \ar[r]  & H^{1}_{C_{\bar{t}}}(C; {\mathbb Z}_C)  \ar[r]  & 0\\
}
$$
and 
$$H^{l}_{C_{\bar{t}}}(C; {\mathbb Z}_C)\simeq  H^{l-1}(C\setminus C_{\bar{t}} ; {\mathbb Z}_C)$$
for $2\leq l\leq r.$
So  by exactness and Fact \ref{fact bf cells} (2) in this case we also have 
\begin{equation*}
H_{C_{\bar{t}}}^{l}(C ; {\mathbb Z}_C) =
\begin{cases} 
{\mathbb Z} \qquad \textmd{if} \qquad l=r\\
\\
\,\,\,\,\,\,\,\,\,\,\,\,\,\,\,\,\,\,\,\,\,\,\,\,\,\,\,\,\,\,\,
\\
0\qquad \textmd{if} \qquad l\neq r.
\end{cases}
\end{equation*}

Now by Theorem \ref{thm Uct} (5), for every definably compact definable subset $A$ of $C$ there is $\bar{t}$ such that $A\subseteq C_{\bar{t}}.$ 
Therefore,  by definition of o-minimal cohomology with definably compact supports, we have:
\begin{displaymath}
\begin{array}{ccl}
H_c^{*}(C ; {\mathbb Z}_C) & = & \lind {A\in c}H^*_A(C; {\mathbb Z}_C)\\
                                &= & \lind {\bar{t}}H^*_{C_{\bar{t}}}(C; {\mathbb Z}_C).\\
\end{array}
\end{displaymath}
and so the result  follows. 

Let $u\in C.$ Since $H_{C_{t_1, \ldots , t_m}}^*(C;{\mathbb Z}_C)\simeq H_c^*(C;{\mathbb Z}_C)$ for every $t_1, \ldots , t_m,$ taking $c_i<t_i$ for all $i$ we have $H_{\{p_C\}}^*(C;{\mathbb Z}_C)\simeq H_c^*(C;{\mathbb Z}_C)$ (using the notation of Theorem \ref{thm Uct} (3)). So it is enough to show that  $H_{\{p_C\}}^*(C;{\mathbb Z}_C)\simeq H_{\{u\}}^*(C;{\mathbb Z}_C)$  

Since $C\subseteq \Pi _{i=1}^m(-_i\frac{b_i}{4},\frac{b_i}{4})$ is a cell of  dimension $r,$ it  is definably homeomorphic, via the restriction of a project of $\Pi _{i=1}^m(-_i\frac{b_i}{4},\frac{b_i}{4})$ to some $r$ coordinates, to an open cell (\cite[Chapter 3, (2.7)]{vdd}). Thus we may assume that $r=m$ and $C$ is open in $\Pi _{i=1}^m(-_i\frac{b_i}{4},\frac{b_i}{4}).$  Recall that, by Fact \ref{fact grp-int grp by 4},   the group-interval operations $x+_iy$, $x-_iy$ and $\frac{x}{2}$ are all defined in each coordinate of 
$ \Pi _{i=1}^m(-_i\frac{b_i}{4},\frac{b_i}{4}).$ So if $u=\langle u_1,\ldots , u_m\rangle $ and $p_C=\langle p_1, \ldots , p_m\rangle ,$ then  $v=\langle -_1p_1+_1u_1, \ldots , -_mp_m+_mu_m\rangle \in \bJ$ is defined. Therefore, by Remark \ref{nrmk grp-int}, there is an open definable neighborhood $W$ of $p_C$ in $C$ such that  right  translation (in $\bJ$) by $v,$ $r_v: W\to r_v(W)\subseteq C,$ sends $p_C$ to $u$ and is a definable homeomorphism. Hence, $r_v^*:H_{\{u\}}^*(r_v(W);{\mathbb Z}_C)\to H_{\{p_C\}}^*(W;{\mathbb Z}_C)$ is an isomorphism. Since by excision, the inclusions $W\to C$ and $r_v(W)\to C$ induce isomorphisms $H_{\{p_C\}}^*(C;{\mathbb Z}_C)\to H_{\{p_C\}}^*(W;{\mathbb Z}_C)$ and $H_{\{u\}}^*(C;{\mathbb Z}_C)\to H_{\{u\}}^*(r_v(W);{\mathbb Z}_C),$ the result follows.
\qed \\

By  Theorem \ref{thm open cells} and  Proposition \ref{prop c-cohom of j-cell} we have:

\begin{cor}\label{cor has orient def jman}
Let $X$ be a definable manifold  of dimension $n$ with definable $\bJ$-bounded charts. Then every   open definable subset $U\subseteq X$ is  a finite union of   open definable subsets $U_1,\dots , U_l \subseteq U$  each of which is definably homeomorphic to a $\bJ$-bounded cell of  dimension  $n$ and such that, for each $i$, we have:
\begin{equation*}
H_c^p(U_i;{\mathbb Z}_X)=
\begin{cases}
{\mathbb Z} \qquad \textmd{if} \qquad p=n\\
\\
\,\,\,\,\,\,\,\,\,\,\,\,\,\,\,\,\,\,\,\,\,\,\,\,\,\,\,\,\,\,\,
\\
0\qquad \textmd{if} \qquad p\neq n. \\
\end{cases}
\end{equation*}
\end{cor}
$\,$\\

By Corollary \ref{cor has orient def jman}  and Definition \ref{defn orient sheaf}  we then have:\\

\begin{cor}\label{cor orient sheaf}
Let $X$ be a definable manifold  of dimension $n$ with definable $\bJ$-bounded charts such that $c$ is a family of definably normal supports. Then $X$ has  orientation sheaf ${\mathcal Or}_X$ which is locally constant with sections given by 
$$\Gamma (U;{\mathcal Or}_X)\simeq \Ho (H_c^n(U;\ZZ_X), \ZZ)$$
for each open definable subset $U\subseteq X.$ Moreover,  for each $x\in X$ we have 
$$\Or _{X,x}\simeq \Ho (H_{\{x\}}^n(X;\ZZ_X), \ZZ).$$
\end{cor}

\pf
Let $x\in X.$ If $V\subseteq X$ is an open definable subset of $X$ such that $x\in V,$ then by Theorem \ref{thm open cells}, there is an open definable subset $U$ of $V$ such that $x\in U$ and $U$ is definably homeomorphic to an open $\bJ$-cell. Therefore 
\begin{align*}
\Or _{X,x} &=\lind {x\in W}\Gamma (W; \Or _X)\\
\end{align*}
where the limit is now over open definable subsets $W$ of $X$ such that $x\in W$ and $W$ is definably homeomorphic to an open $\bJ$-cell. 

For such $W$'s we also have
\begin{align*}
\Gamma (W; \Or _X)&\simeq  \Ho (H_{c}^n(W;\ZZ_X), \ZZ)\\
& \simeq \Ho (H_{\{x\}}^n(W;\ZZ_X), \ZZ)\\
& \simeq  \Ho(H_{\{x\}}^n(X;\ZZ_X), \ZZ).\\
\end{align*}
by Definition \ref{defn orient sheaf}, Proposition \ref{prop c-cohom of j-cell} and excision isomorphism. Hence we obtain $\Or _{X,x}\simeq \Ho (H_{\{x\}}^n(X;\ZZ_X), \ZZ).$ \qed \\

\end{subsection}

\end{section}

\begin{section}{Applications to  definably compact abelian groups}\label{section  torsions}
In this section we prove our results about definably compact abelian groups as explained in the Introduction.

\begin{subsection}{Definably compact groups and their intrinsic o-minimal fundamental groups}\label{subsection def grps and gr-int}
The following result establishes the connection between definable groups and cartesian products of definable groups-intervals:

\begin{fact}\cite[Theorem 3]{epr}\label{fact grp epr}
If $G$ is a definable group, then there is a definable injection $G\to \Pi _{i=1}^mJ_i$, where each $J_i\subseteq M$ is a definable group-interval.
\end{fact}

In particular, by \cite{p1} and \cite[Corollary 2.3]{et} (definable normality of definable groups) and also \cite[Theorem 2.1]{ps}, we have that if $G$ is a definably compact definable group, then there is a  cartesian product  $\bJ=\Pi _{i=1}^mJ_i$ of definable group-intervals such that $G$ is a definable manifold   with definable $\bJ$-bounded charts such that $c$ is a family of definably normal supports. \\

Below we fix such  $\bJ=\Pi _{i=1}^mJ_i$ and since $\bJ$ is constructed from $G$  we will call $\pi _1^{\bJ}(G)$ the {\it intrinsic o-minimal fundamental group of $G$}. \\

By Remark \ref{nrmk covers main}  we can extend to intrinsic o-minimal fundamental groups of definably compact definable groups in ${\mathbb M}$ (an arbitrary o-minimal structure), the results about o-minimal fundamental groups already proved in o-minimal expansions of ordered fields (\cite[Section 2]{eo}) or in o-minimal expansions of ordered groups (\cite{EdPa}).\\

\begin{thm}\label{thm uni cover of g}
Let  $G$ be a definably compact, definably connected definable group. Then there exists a universal locally definable covering homomorphism $\tilde{p}:\tilde{G}\into G$ where $\tilde{G}$  is a locally definable manifold with definable $\bJ$-bounded charts. Moreover, the intrinsic o-minimal fundamental group  $\pi _1^{\bJ}(G)$ of $G$ is  abelian  and finitely generated.\\
\end{thm}

\pf
By Remark \ref{nrmk covers main} and \cite[Theorem 1.2]{eep} there is a universal locally definable covering map $\tilde{p}:\tilde{G}\to G$ in the category of locally definable manifolds with definable $\bJ$-bounded charts. So  $\tilde{G}$  is a locally definable manifold with definable $\bJ$-bounded charts. By \cite[Proposition 2.28]{eep} the definable group operations of $G$ can be lifted to locally definable group operations on $\tilde{G}$ making  $\tilde{G}$ a locally definable group and $\tilde{p}:\tilde{G}\to G$ a locally definable homomorphism (compare with the proof of \cite[Claims 3.9 and 3.10]{EdPa}).

By Remark \ref{nrmk covers main} and \cite[Theorem 1.1]{eep} the o-minimal $\bJ$-fundamental group  $\pi _1^{\bJ}(G)$ of $G$ is   finitely generated. As in  \cite[Lemma 2.3]{eo}, $\pi _1^{\bJ}(G)$ is abelian.
\qed \\

The following is proved in exactly the same way as the proof of its analogue in o-minimal expansions of ordered fields (\cite[Theorem 2.1]{eo}): 

\begin{thm}\label{thm pi g and tor}
Let $G$ be a definably compact, definably connected definable abelian group. Then there is $s\in {\NN}$ such that:
\begin{enumerate}
\item[(a)]
the intrinsic o-minimal fundamental group  $\pi _1^{\bJ}(G)$ of $G$ is isomorphic to ${\ZZ}^s,$ and 
 \item[(b)]
the subgroup $G[k]$ of $k$-torsion points of $G$ is isomorphic to $({\ZZ}/k{\ZZ})^s$, for each $k\in {\NN}.$
 \end{enumerate}
\end{thm}

\pf
Since \cite[Proposition 2.10]{eo} holds in arbitrary o-minimal structures, $p_k:G\to G:x\mapsto kx$ is a definable covering homomorphism with $G[k]\simeq {\rm Aut}(p_k:G\to G)$ (the group of definable homeomorphisms of  coverings). By 
\cite[Theorem 3.9]{eep} (and Remark \ref{nrmk covers main}) we have also $G[k]\simeq \pi _1^{\bJ}(G)/p_{k*}(\pi _1^{\bJ}(G)).$ Since \cite[Lemmas 2.3 and 2.4]{eo} hold here as well, $p_{k*}(\pi _1^{\bJ}(G))=k\pi _1^{\bJ}(G).$ Finally, by 
\cite[Corollary 2.17]{eep} (and Remark \ref{nrmk covers main}), $p_{k*}:\pi _1^{\bJ}(G) \to \pi _1^{\bJ}(G)$ is injective and so the finitely generated abelian group $\pi _1^{\bJ}(G)$ is free.
\qed \\

Finally, by Remark \ref{nrmk covers main} and \cite[Theorem 4.16]{eep} we also have the Hurewicz theorem:

\begin{thm}\label{thm hure thm}
Let $G$ be a definably compact, definably connected definable abelian group. Then 
$$ \Ho (\pi _1^{\bJ}(G)^{\rm op}, \ZZ) \simeq \check{H }^{1}(G;\ZZ).$$\\
\end{thm}

\end{subsection}

\begin{subsection}{Orientability of definably compact groups}\label{subsection orient def grp}
Here we show that definably compact definable groups are orientable.\\

\begin{thm}\label{thm def comp grp is orient}
Let $G$ be a definably connected, definably compact, definable group of positive dimension $n$. Then $G$ is orientable. In particular, $H^n(G;\ZZ_G)\simeq \ZZ.$
\end{thm}

\pf
By Corollary \ref{cor orient sheaf}, $G$ has  orientation sheaf ${\mathcal Or}_G$ which is locally constant with sections given by 
$$\Gamma (U;{\mathcal Or}_G)\simeq \Ho (H_c^n(U;\ZZ_G), \ZZ)$$
for each open definable subset $U\subseteq G.$ Moreover,  for each $x\in G$ we have 
$$\Or _{G,x}\simeq \Ho (H_{\{x\}}^n(G;\ZZ_G), \ZZ).$$

By Fact \ref{fact or and lief cov}, there is a canonically associated locally definable $\ZZ$-covering map
$$w_{\oO}:W_{\oO}\to G$$
where $W_{\oO}=\bigsqcup _{x\in G}{\Or}_{G,x}$ and $w_{\oO}(s_x)=x,$ such that $G$ is orientable if and only if the locally definable $\ZZ$-covering map $w_{\oO}:W_{\oO}\to G$ is trivial. Moreover, if $U_1,\ldots , U_l$ is an admissible cover of $G$ by open definable subsets such that for each $j=1,\ldots, l,$ the restriction ${\Or}_{G|U_j}\simeq {\mathbb Z}_{G|U_j},$ then $w_{\oO}:W_{\oO}\to G$ is a locally definable $\ZZ$-covering map trivial over ${\mathcal U}=\{U_j\}_{j=1}^l$ with, for each $j=1,\ldots, l:$
\begin{itemize}
\item
$w^{-1}_{\oO}(U_j)=\bigsqcup _{m\in {\ZZ}}U_j^m;$
\item
each $w_{\oO|U_j}:U_j^m\to U_j$ is a definable homeomorphism
\end{itemize}
where for $m\in {\ZZ}$, $U_j^m=\{s_x:\Gamma (U_j; {\Or}_G)\simeq \ZZ:s\mapsto m,\,\,\,x\in U_j\}.$

Let $e_G\in G$ be  the identity element and suppose that  $e_G\in U_1.$ Let $\omega _{e_G}\in  \Or _{G,x}$ corresponding to $1\in \ZZ$ under the isomorphism  $\Or _{G,x}\simeq \ZZ $ induced by the isomorphism ${\Or}_{G|U_1}\simeq {\mathbb Z}_{G|U_1}.$

For $z\in G$ let $L_z:G\to G:u\mapsto zu$ be the left translation by $z.$  Then $L_z^*:H_{\{z\}}^n(G;\ZZ _G)\to H_{\{e_G\}}^n(G;\ZZ _G)$ is an isomorphism. So we can define $\omega _z\in \Or _{G,z}$ to be the image $(L_z^*)^{\, \vee }(\omega  _{e_G})$ of $\omega  _{e_G}$ under the dual of the isomorphism $L_z^*.$ 

\begin{clm}\label{clm gen}
For each $j$, there is a unique generator 
$$\omega _j\in \Gamma (U_j;{\mathcal Or}_G)\simeq \Ho (H_c^n(U_j;\ZZ_G), \ZZ)\simeq \ZZ$$ 
such that  $\omega _z=(\omega _j)_{z}$ for all $z\in U_j.$ 
\end{clm}

Since $U_j$ is definably homeomorphic to an open $\bJ$-cell $C_j,$ if  $u_j\in U_j$ then by excision and Proposition \ref{prop c-cohom of j-cell} we have
$$H_{\{u_j\}}^n(G;\ZZ _G)\simeq H_{\{u_j\}}^n(U_j;\ZZ _G)\simeq H_c^n(U_j;\ZZ _G)\simeq \ZZ.$$ 
Therefore, if  $z\in U_j$, then we have a commutative  diagram of isomorphisms (where the dashed arrow is introduced to make the diagram commutative):
$$
\xymatrix{
 H^n_{u_j}(G;\ZZ _G)  \ar[r]^{\simeq } \ar[d] ^{L_{u_jz^{-1}}^*}  &  H^n_{u_j}(U_j;\ZZ_G) \ar[r]^{\simeq}  &   H_c^n(U_j; \ZZ_G)\ar[d]^{{\rm id}} \\
 H^n_{z}(G;\ZZ _G) \ar[r]^{\simeq } & H^n_{z}(U_j;\ZZ _G) \ar@{-->}[r] & H_c^n(U_j;\ZZ_G).
}
$$
Since $L_{u_j}=L_{u_jz^{-1}}\circ L_z$ we have $(L_{u_j}^*)^{\, \vee }=(L_{u_jz^{-1}}^*)^{\, \vee}\circ (L_z^*)^{\, \vee }$ and so $\omega _{u_j}=(L_{u_jz^{-1}}^*)^{\, \vee}(\omega _z).$ This together with the (dual) of the commutative diagram above shows that there is a unique generator $\omega _j\in \Gamma (U_j;{\mathcal Or}_G)\simeq \Ho (H_c^n(U_j;\ZZ_G), \ZZ)\simeq \ZZ$ such that $\omega _{u_j}=(\omega _j)_{u_j}$ and $\omega _z=(\omega _j)_{z}.$ 
\qed  \\

\begin{clm}
Let $s:G\to W_{\oO}$  be the map given by $s(z)=\omega  _z.$ Then $s:G\to W_{\oO}$ is a continuous locally definable section to $w_{\oO}:W_{\oO}\to G$ (i.e. such that $w_{\oO}\circ s={\rm id}_G$).
\end{clm}

We have $w_{\oO}\circ s={\rm id}_G$ and moreover, since $G$ is definably connected, by Claim \ref{clm gen}, for each $j,$ we have $s(U_j)=U_j^1=\{s_x:\Gamma (U_j; {\Or}_G)\simeq \ZZ:s\mapsto 1,\,\,\,x\in U_j\}.$ Hence, $s_{|U_j}:U_j\to U_j^1$ is the inverse to the definable homeomorphism $w_{\oO|U_j}:U_j^1\to U_j.$ Therefore, $s$ is continuous and locally definable as required. \qed \\

Finally, if $\alpha :W_{\oO}\times {\ZZ} \to W_{\oO}$ is the locally definable ${\ZZ}$-action making $w_{\oO}:W_{\oO}\to G$ into a locally definable ${\ZZ}$-covering map, then $W_{\oO}\to G\times {\ZZ}: v\mapsto (w_{\oO}(v), l_v)$ where $l_v\in \ZZ$ is the unique element such that $v=\alpha (s(w_{\oO}(v)), l_v)$ is a locally definable homeomorphism of locally definable ${\ZZ}$-coverings  showing that $w_{\oO}:W_{\oO}\to G$  is trivial. \qed \\


\end{subsection}

\begin{subsection}{The Hopf algebra of a definably compact group}\label{subsection hopf def grp}
Here we show that the o-minimal cohomology $H^*(G;k_G)$ of a definably connected, definably compact definable group $G$ with coefficients in a field $k$  is a connected, bounded, Hopf algebra over $k$ of finite type.\\

For the rest of this subsection let $k$ be a field. We start with a couple of  observations required below:\\

By  the isomorphism $\mod (k_{X_{\df}})\simeq \mod (k_{\tilde{X}})$ (Theorem \ref{thm main iso on sheaves}) and \cite[Chapter II, Section 7 and (8.2) ]{b} 
we have:

\begin{fact}\label{fact cup}
Let $X$ be a definable space. Then there is a cup product operation
$$\cup : H^p(X;k_X)\otimes H^q(X;k_X)\to H^{p+q}(X;k_X)$$
making $H^*(X;k_X)$ into a graded, associative, skew-commutative $k$-algebra with unit in $H^0(X;k_X)$. This product is functorial and the algebra is connected if $X$ is definably connected. \\
\end{fact}



\begin{fact}[K\"unneth formula]\label{fact kunn grps}
Let $G$ be a definably compact definable group. Then there is a natural isomorphism
$$H^*(G\times G;k_{G\times G})\simeq \bigoplus _{p+q=*}(H^p(G;k_G)\otimes H^q(G;k_G)).$$ 
\end{fact}

\pf
Consider the full subcategory of locally closed definable subsets of definably compact definable groups. We claim that this category satisfies conditions (A0), (A1) and (A2) and definably compact groups satisfy condition (A3) listed  in the Introduction.  As explained in the Introduction, after the work developed in \cite{ep3}, this ensures that we have K\"unneth formula for $G\times G.$

(A0) follows from that fact that  a product of  locally closed definable subsets of a cartesian product of  definably compact definable groups  is also a locally closed definable subset of a  definably compact definable group; (A1) follows from Corollary \ref{cor def jman basis norm}; (A2) follows since: (i)  a  definably compact group is  definably normal (\cite[Corollary 2.3]{et}) and (ii) a locally closed definable subset of a cartesian product of a given definably compact definable group  has a definably normal completion, namely its closure; (A3) was proved in \cite[Theorem 1.1]{ep2}. 
\qed \\

\begin{nrmk}\label{nrmk kunn grps funct}
{\em
Since the  K\"unneth isomorphism (Fact \ref{fact kunn grps}) is natural, we have: 
\begin{itemize}
\item
If $f: G\to G$ and $g:H\to H$ are definable continuous maps between definably compact definable groups, then  we can make the identification  $(f\times g)^*=f^*\otimes g^*.$ 
\item
If $q_i:G\times G \to G$ ($i=1,2$) is the projection onto the $i$-th coordinate,  then using $\{0\}\times G=G=G\times \{0\},$ it follows that for $z\in H^p(G;k_G)$  we can make the identifications  $q_1^*(z)=z\otimes 1$ and $q_2^*(z)=1\otimes z.$
 \end{itemize}
}
\end{nrmk}

We can now prove the main result of this subsection:\\

\begin{thm}\label{thm hopf grps}
Let $G$ be a definably connected, definably compact definable group.  Then the o-minimal sheaf cohomology $H^*(G;k_G)$ of $G$ with coefficients in  $k$ is a connected, bounded, Hopf algebra over $k$ of finite type. Moreover, if ${\rm char}(k)=0$, then we have a Hopf algebra isomorphism
$$H^*(G;k_G)\simeq \bigwedge [y_1, \ldots , y_r]_{k}$$
with the exterior algebra  with the $y_i$'s of odd degree and primitive. 
\end{thm}

\pf
By Fact \ref{fact cup}, $H^*(G;k_G)$ is a connected, graded, associative, skew-commutative $k$-algebra with unit in $H^0(G;k_G)$.

Let  $m:G\times G\to G$ the multiplication map. Let
$$\mu :H^*(G;k_G)\to  \bigoplus _{p+q=*}(H^p(G;k_G)\otimes H^q(G;k_G))$$
be the composition of $m^*:H^*(G;k_G)\to  H^*(G\times G;k_{G\times G})$ with the isomorphism $H^*(G\times G;k_{G\times G})\simeq \bigoplus _{p+q=*}(H^p(G;k_G)\otimes H^q(G;k_G))$ given by the K\"unneth formula (Fact \ref{fact kunn grps}). From the properties of $m$, it  is standard to show that $\mu $ is a co-multiplication making $H^*(G;k_G)$ into a Hopf algebra over $k$. Compare with \cite[Corollary 3.5]{eo}. 

Since $G$ is definably normal (\cite[Corollary 2.3]{et}), by \cite[Proposition 4.2]{ejp1} we have  $H^p(G;k_G)=0$ for all $p>\dim G$ and therefore $H^*(G;k_G)$ is a bounded Hopf algebra. By \cite[Theorem 1.2]{ep2}, it is of finite type.

In the case ${\rm char}(k)=0,$ the description of the Hopf algebra $H^*(G;k_G)$ follows from the Hopf-Leray theorem. See \cite[Corollary 3.6]{eo} for details.
\qed \\

\end{subsection}

\begin{subsection}{Computing the torsion subgroups}\label{subsection torsion subgrps}
Here we  compute the torsion subgroups of a definably compact abelian definable group.\\

Below we will omit  the subscript on the field ${\mathbb Q}$ when we consider the constant sheaf it determines on a definably compact definable group $G$. We also consider $G$ with a fixed orientation (Theorem \ref{thm def comp grp is orient}). 

\begin{lem}\label{pkdeg} 
Let $G$ be a definably compact, definably connected, abelian definable group. For each $k>0$, consider the map $p_k\colon G\to G\colon x\mapsto kx$. Then we have  ${\rm deg}\,p_k\leq |p_k^{-1}(0)|,$ where $0$ is the neutral element of $G$.
\end{lem}

\pf
 Fix  a $k>0$. By Corollary \ref{pointdeg}, ${\rm deg}\,p_k={\rm deg}_0p_k$. Also, by \cite[Corollary 2.12]{eo} (which holds in arbitrary o-minimal structures), we know that  the homomophism  $p_k\colon G\to G$ is a definable covering map. Let $G=\bigcup_{l\in L}U_l$ be  as in the definition of definable covering map. Fix $l_o\in L$ such that $0\in U_{l_o}$ and let $Y=p_k^{-1}(U_{l_o})$. Now consider the map
$f=(p_k)_{|Y}:Y\to G$. Note that, by the dual of the excision, ${\rm deg}_0p_k={\rm deg}_0f$. On the
other hand,  if  we let $Y=\sqcup_{\lambda=1}^m Y_{\lambda}$ with the $Y_{\lambda}$'s being the definably connected components of $Y$ so that $f^{\lambda}=f_{|Y_{\lambda}}:Y_{\lambda }\to G$ is a homeomorphism onto an open subset of $G$, namely $U_{l_o}$. Now the data $f:Y\to G$, $\{ 0\}\subset G$ and $Y=\sqcup_{\lambda=1}^m Y_{\lambda}$ satisfy the hypothesis of Lemma \ref{sumdeg}, therefore  we can conclude that  ${\rm deg}_0f=\sum_{\lambda=1}^m {\rm deg}_0f^{\lambda}$. Finally, by Lemma \ref{homeomdeg}, ${\rm deg}_0f^{\lambda}=\pm 1$, for each $\lambda$, and hence ${\rm deg}_0f\leq m=|f^{-1}(0)|$. By the above, this last relation means  ${\rm deg}\,p_k\leq |p_k^{-1}(0)|$.
\qed \\

By Theorem \ref{thm hopf grps} we have a Hopf algebra isomorphism
$$H^*(G;\QQ)\simeq \bigwedge [y_1, \ldots , y_r]_{\QQ}$$
with the exterior algebra  with the $y_i$'s of odd degree and primitive. This means that $\mu (y_i)=y_i\otimes 1 + 1\otimes y_i$ for each $i=1,\ldots ,r $ where the co-multiplication 
 $$\mu =m^*:H^*(G;\QQ)\to  H^*(G\times G; \QQ)\simeq \bigoplus _{p+q=*}(H^p(G;\QQ)\otimes H^q(G;\QQ))$$
is given by the composition of the homomorphism  $m^*:H^*(G;\QQ)\to  H^*(G\times G;\QQ)$ induced by the multiplication map $m:G\times G \to G$ on $G$ with the isomorphism $H^*(G\times G;\QQ)\simeq \bigoplus _{p+q=*}(H^p(G;\QQ)\otimes H^q(G;k_G))$ given by the K\"unneth formula (Fact \ref{fact kunn grps}).
 
We call an element $x\in H^*(G;\QQ)$ a {\it monomial of length $l$} if $x=y_{i_1}\cup \cdots \cup y_{i_l}$ where $1\leq i_1<\ldots <i_l\leq r.$\\

\begin{lem}\label{pk*}
Let $G$ be a definably connected, definably compact, definable  group. For each $k>0$, consider  the definable continuous map $p_k:G\to G:a\mapsto a^k$, for each $a\in G$. Then,
the map $p_k^*: H^*(G;\QQ)\to H^*(G;\QQ)$ sends  each monomial $x$ of length $l$ to $k^lx$.
\end{lem}

\pf 
First we prove by induction on $k$ that, for   $y\in\{ y_1, \ldots, y_r\},$  we have $p_k^*(y)=ky.$ 
 
 For $k=1,$  we have $p_k=\id $ and so this case is trivial.  For the induction step, using $p_{k+1}=m\circ(p_k\times \id)\circ \Delta $ where $\Delta : G\to G\times G$ is the diagonal map in $G$,  we have   
\begin{align*}
p_{k+1}^*(y) & =(m\circ(p_k\times \id)\circ \Delta )^*(y)\\
& =(\Delta ^*\circ(p_k\times \id )^*\circ m^*)(y)\\
&=\Delta ^*\circ(p_k^*\otimes \id)(y\otimes 1+ 1\otimes y)\\
& =\Delta ^*(ky\otimes 1 + 1\otimes y) \\
&= \Delta ^*(q_1^*(ky) + q_2^* (y)) \\
&=ky + y \\
&=(k+1)y. 
\end{align*}
In these equalities we used Remark \ref{nrmk kunn grps funct} and $q_i\circ \Delta =\id$ where $q_i:G\times G\to G$ ($i=1,2$) is the projection onto the $i$th coordinate.

Finally,  we get  $p_{k+1}^*(x)=(k+1)^lx$, for each $k>0$,  since $p_{k+1}^*$ is an algebra morphism.
\qed \\

\medskip
\noindent
{\bf Proof of Theorem \ref{thm main thm}:}
Let $G$ be a definably connected definably compact  abelian group of dimension $n$. Consider also the Hopf algebra $H^*(G;\QQ)\simeq \bigwedge [y_1, \ldots , y_r]_{\QQ}$ of $G.$  

As we saw in Remark \ref{nrmk fund class} the orientation class $\mu _G\in \Gamma (G; \Or _G)$ determines the fundamental class $\zeta _G\in H^n(G; \ZZ _G)^{\, \vee }.$ Let $\zeta _G^{\, \vee }\in H^n(G; \ZZ _G)$ be the dual of $\zeta _G$ and let $\omega _G \in H^n(G;\QQ)$ be the image of $\zeta _G^{\, \vee }$ under the isomorphism 
$$H^n(G; \ZZ _G)\otimes \QQ \simeq H^n(G; \QQ)$$ 
given by the universal coefficients formula (\cite{ep3}).

Now fix a $k>0$, and  consider the definable continuous map $p_k:G\to G: a\mapsto ka$. By definition of degree of a map we obtain $p_k^*(\omega_G)=({\rm deg}\,p_k)\omega_G$. Since $\omega_G$ generates $H^n(G;\QQ)$, and $0\not=y_1\cup \cdots \cup y_r\in H^n(G;\QQ)$ we can suppose $\omega_G=y_1\cup \cdots\cup  y_r$. By Lemma \ref{pk*}, $p_k^*(\omega_G)=k^r\omega_G$, and so ${\rm deg} \,p_k=k^r$.

On the other hand, by Theorem \ref{thm pi g and tor}, there is an $s\geq 0$ such that $\pi _1^{\bJ}(G)\simeq  \ZZ ^s$ and $p_k^{-1}(0)=G[k]\simeq  (\ZZ/k\ZZ)^s$. By Corollary \ref{pkdeg}, ${\rm deg}\,p_k\leq|(p_k)^{-1}(0)|=k^s$, and hence $r\leq s$. By the Hurewicz theorem (Theorem \ref{thm hure thm}),  
$$ \Ho (\pi _1^{\bJ}(G)^{\rm op}, \ZZ) \simeq \check{H }^{1}(G;\ZZ),$$
and on the other hand, since $G$ is definably normal (\cite[Corollary 2.3]{et}),  $\check{H}^1(G;\ZZ)\simeq H^1(G; \ZZ _G)$ (\cite[Proposition 4.1]{ejp1}) and hence 
$$\Ho (\pi _1^{\bJ}(G)^{\rm op}, \ZZ)\otimes \QQ \simeq  H^1(G;\QQ).$$ 
Therefore, since $H^1(G;\QQ)$ is a subspace of $H^*(G;\QQ)$ (and the elements of $H^1(G;\QQ)$ cannot be decomposable), among $\{ y_1,\dots,y_r\}$ there must be exactly $s$ elements of degree one.  Hence $s=r$ and all $y_i$'s are  of degree one. Finally, since $\omega _G=y_1\cup \cdots \cup y_r\in H^n(G;\QQ)$ we must have $s=r=n$. 
 Therefore, 
\begin{enumerate}
\item[(a)] 
$\pi _1^{\bJ}(G)\simeq {\ZZ}^n$;

\item[(b)]  
$G[k]\simeq({\ZZ}/k{\ZZ})^n$, and 
 
\item[(c)]   
$H^*(G;\QQ)\simeq \bigwedge [y_1, \ldots , y_n]_{\QQ}$ with the $n$ generators of degree one.
\end{enumerate}
 \qed \\

\end{subsection}

\end{section}

\end{document}